\ifx\shlhetal\undefinedcontrolsequence\let\shlhetal\relax\fi
\documentclass{amsart}
\usepackage{amsmath} 
\usepackage{amssymb}

\newtheorem{theorem}{Theorem}[section] 
\newtheorem{claim}{Claim}[theorem]
\newtheorem{lemma}[theorem]{Lemma} 
\newtheorem{proposition}[theorem]{Proposition} 
\newtheorem{corollary}[theorem]{Corollary} 

\theoremstyle{definition}
\newtheorem{definition}[theorem]{Definition}
\newtheorem{problem}{Problem}[section]

\theoremstyle{remark}
\newtheorem{remark}[theorem]{Remark}
\newtheorem{conclusion}[theorem]{Conclusion}
\newtheorem{hypothesis}[theorem]{Hypothesis}

\newcommand{\bB}{{\mathbb B}}

\newcommand{\bP}{{\mathbb P}}
\newcommand{\cA}{{\mathcal A}}
\newcommand{\cB}{{\mathcal B}}
\newcommand{\cD}{{\mathcal D}}

\newcommand{\cU}{{\mathcal U}}
\newcommand{\hd}{{\rm hd}}
\newcommand{\hL}{{\rm hL}}
\newcommand{\Ult}{{\rm Ult}}

\newcommand{\comp}{\circ}
\newcommand{\forces}{\Vdash}
\newcommand{\V}{{\bf V}}
\newcommand{\rest}{\restriction}
\newcommand{\cf}{{\rm cf}}

\newcommand{\otp}{{\rm otp}}
\newcommand{\cl}{{\rm cl}}
\newcommand{\id}{{\rm id}}
\newcommand{\xs}{{\mathcal X}_S}
\newcommand{\qs}{{\mathbb Q}_S}
\newcommand{\Ps}{{\mathbb P}_S}
\newcommand{\vare}{\varepsilon}
\newcommand{\dcD}{{\dot{\cD}}}
\newcommand{\cof}{{\rm cof}}
\newcommand{\lex}{{\rm lex}}
\newcommand{\lh}{{\rm lh}}

\title{Forcing for $\hL$ and $\hd$}
\author{Andrzej Ros{\l}anowski}
\address{Department of Mathematics\\
 University of Nebraska at Omaha\\
 Omaha, NE 68182-0243, USA\\
 and Mathematical Institute of Wroclaw University\\
 50384 Wroclaw, Poland} 
\email{roslanowski@unomaha.edu}
\urladdr{http://www.unomaha.edu/$\sim$aroslano}
\thanks{The first author thanks the KBN (Polish Committee of  Scientific
 Research) for partial support through grant 2 P03 A 01109.} 

\author{Saharon Shelah}
\address{Institute of Mathematics\\
 The Hebrew University of Jerusalem\\
 91904 Jerusalem, Israel\\
 and  Department of Mathematics\\
 Rutgers University\\
 New Brunswick, NJ 08854, USA}
\email{shelah@math.huji.ac.il}
\urladdr{http://www.math.rutgers.edu/$\sim$shelah}
\thanks{The research of the second author was partially supported by the
 Israel Science Foundation. Publication 651} 


\begin{document}

\begin{abstract}
The present paper addresses the problem of attainment of the supremums
in various equivalent definitions of hereditary density $\hd$ and
hereditary Lindel\"of degree $\hL$ of Boolean algebras.  We partially answer
two problems of J.~Donald Monk, \cite[Problems 50, 54]{M2}, showing
consistency of different attainment behaviour and proving that (for the
considered variants) this is the best result we can expect.
\end{abstract}

\maketitle

\section{Introduction}
We deal with the attainment problem in various definitions of two cardinal 
functions on Boolean algebras: the hereditary density $\hd$ and the
hereditary Lindel\"of degree $\hL$. These two cardinal functions are closely 
related, as it is transparent when we pick the right variants of (equivalent)
definitions. Also they both are somewhat related to the spread $s$ of
Boolean algebras. So, for a Boolean algebra $\bB$, we define
\begin{itemize}
\item $s(\bB)=\sup\{\kappa:$ there is an ideal--independent sequence of
length $\kappa\;\}$,
\item $\hd(\bB)=\sup\{\kappa:$ there is a left--separated sequence of length
$\kappa\;\}$,
\item $\hL(\bB)=\sup\{\kappa:$ there is a right--separated sequence of length
$\kappa\;\}$.
\end{itemize}
Let us recall that a sequence $\langle a_\xi:\xi<\kappa\rangle$ of elements
of a Boolean algebra is 
\begin{itemize}
\item {\em ideal--independent\/} if $a_\xi\nleq\bigvee\limits_{\zeta\in w}
a_\zeta$ for each $\xi<\kappa$ and a finite set $w\subseteq \kappa\setminus
\{\xi\}$,   
\item {\em left--separated\/} if $a_\xi\nleq\bigvee\limits_{\zeta\in w}
a_\zeta$ for each $\xi<\kappa$ and a finite set $w\subseteq \kappa\setminus
(\xi+1)$,   
\item {\em right--separated\/} if $a_\xi\nleq\bigvee\limits_{\zeta\in w}
a_\zeta$ for each $\xi<\kappa$ and a finite set $w\subseteq\xi$.   
\end{itemize}
The above definitions of the three cardinal functions are of special use,
see e.g. \cite[\S 1]{RoSh:534}. However, neither these definitions explain
the names of the functions, nor they are good enough justifications for the 
interest in them. But all three functions originate in the cardinal
functions of the topological space $\Ult(\bB)$ (of ultrafilters on $\bB$). 
And thus, for a Boolean algebra $\bB$, we may define (or prove that the
following equalities hold true):  
\begin{itemize}
\item $s(\bB)=\sup\{|X|: X\subseteq\Ult(\bB)$ is discrete in the relative
topology $\}$,
\item $\hd(\bB)=\sup\{d(X): X\subseteq\Ult(\bB)\}$, where\\
$d(X)=\min\{|Y|:Y\subseteq X$ is dense in $X\;\}$,
\item $\hd(\bB)=\sup\{L(X): X\subseteq\Ult(\bB)\}$, where\\
$L(X)=\min\{\kappa:$ every open cover of $X$ has a subcover of size $\leq
\kappa\;\}$.
\end{itemize}
The respective pairs of cardinal numbers are defined using $\sup$, so even
if we know that they are equal we still may expect different attainment
properties: one of the families of cardinals may have the largest member
while the other not. Also we may ask if the $\sup$ has to be
attained. Situation may seem even more complicated if one notices that there
are more than just two equivalent definitions of the cardinal functions
$s,\hd,\hL$: Monk \cite{M2} lists six equivalent definitions for spread (see
\cite[Theorem 13.1]{M2}), nine definitions for $\hd$, and nine for $\hL$
(see \cite[Theorems 16.1, 15.1]{M2}). Fortunately, there is a number of
dependencies here. 

First, all of the equivalents of spread have the same attainment
properties. Moreover, spread is always attained for singular strong limit
cardinals and for singular cardinals of countable cofinality (for these and
related results see Hajnal and Juh\'asz \cite{HaJu67}, \cite{HaJu69b},
\cite{HaJu69a}, Juh\'asz \cite{Ju71}, \cite{Ju80}, Roitman \cite{Ro75},
Kunen and Roitman \cite{KuRo77}, Juh\'{a}sz and Shelah \cite{JuSh:231}). 
Then Shelah \cite{Sh:233} proved that $2^{\cf(s(\bB))}<s(\bB)$ implies that
the spread is attained (see \ref{233ons} here). Finally, it is shown in
Shelah \cite[\S 4]{Sh:641} that, e.g., if $\mu$ is a singular strong limit
cardinal such that $\mu<\cf(\lambda)<\lambda\leq 2^\mu$, then there is a
Boolean algebra $\bB$ such that $|\bB|=s(\bB)=\lambda$ and the spread is not
obtained. Thus, to some extend, the problem of attainment for spread is
settled. 

Many of the results mentioned above can be carried out for (some) variants
of $\hd$ and $\hL$. However, the difference between these two cases and the
case of the spread is that the various equivalent definitions of the
respective cardinal function might have different attainment properties.

Let us introduce some of the equivalents of $\hL$, $\hd$. They will be
called $\hL_{(n)}$, $\hd_{(n)}$, with the integer $n$ referring to the
respective cardinal $\kappa_n$ as used in the proofs of \cite[15.1 and
16.1]{M2}, respectively. Also, we will have $\hd^+_{(n)}$ and $\hL_{(n)}^+$
to have proper language to deal with the attainment questions. Let us start
with the hereditary Lindel\"of degree $\hL$. First, for a topological space
$X$ we define the Lindel\"of degree ${\rm L}(X)$ of the space $X$ as
\[{\rm L}(X)=\min\{\lambda:\mbox{every open cover of $X$ has a subcover of
size }\leq\lambda\;\}.\]
\begin{definition}
\label{hLbasic}
Let $\bB$ be an infinite Boolean algebra. For an ideal $I$ in a Boolean
algebra $\bB$ we let 
\[\cof(I)=\min\{|A|:A\subseteq I\mbox{ and }(\forall b\in I)(\exists a\in A) 
(b\leq a)\}.\]
Now we define
\[\begin{array}{l}
\hL_{(0)}^{(+)}(\bB)=\sup\{{\rm L}(X)^{(+)}: X\mbox{ is a subspace of }
\Ult(\bB)\;\},\\
\hL_{(1)}^{(+)}(\bB)=\sup\{\cof(I)^{(+)}: I\mbox{ is an ideal of }\bB\;\},\\
\hL_{(7)}^{(+)}(\bB)=\sup\{\kappa^{(+)}:\mbox{ there is a right--separated
sequence }\langle a_\xi:\xi<\kappa\rangle\mbox{ in }\bB\;\}.
  \end{array}\]
The superscript ``$(+)$'' in the above definitions means that each of the
formulas has two versions: one with ``$+$'' and one without it. 
\end{definition}
The cardinals mentioned in \ref{hLbasic} are among those listed in
\cite[Theorem 15.1]{M2}, and so $\hL_{(0)}(\bB)=\hL_{(1)}(\bB)=\hL_{(8)}(
\bB)$. The attainment properties can be described using the versions with
``$+$'': $\hL^+_{(i)}(\bB)=\hL_{(i)}(\bB)$ means that the supremum is not
obtained; $\hL^+_{(i)}(\bB)=\hL^+_{(j)}(\bB)$ means that the respective two
definitions of $\hL$ have the same attainment for $\bB$. It is not difficult
to note that 
\[\hL^+_{(7)}(\bB)=\hL_{(7)}(\bB)\quad\Rightarrow\quad\hL^+_{(1)}(\bB)=
\hL_{(1)}(\bB)\]
and 
\[\hL_{(0)}(\bB)=\hL_{(0)}^+(\bB)\mbox{ is a regular cardinal }\quad
\Rightarrow\quad \hL_{(7)}^+(\bB)=\hL_{(7)}(\bB)\]
(and the attainment of $\hL$ in senses not listed in \ref{hLbasic} can be
reduced to those three; see \cite[p. 190, 191]{M2} for details). Also, if
$\hL(\bB)$ is a strong limit cardinal or if it has countable cofinality,
then $\hL_{(7)}(\bB)<\hL^+_{(7)}(\bB)$ (see Juh\'asz \cite[4.2, 4.3]{Ju80}). 

In \ref{233forhL} we will show that if $\hL(\bB)$ is a singular cardinal
such that $2^{\cf(\hL(\bB))}<\hL(\bB)$, then $\hL_{(0)}^+(\bB)=\hL_{(1)}^+
(\bB)=\hL_{(7)}^+(\bB)=(\hL(\bB))^+$. Thus, e.g., under GCH, the sups in all 
equivalent definitions of $\hL$ are attained at singular cardinals. Next, in
section 3, we use forcing to show that, consistently, there is a Boolean
algebra $\bB$ such that 
\[\hL_{(7)}(\bB)<\hL^+_{(7)}(\bB)\quad\mbox{ and }\quad\hL^+_{(1)}(\bB)=
\hL_{(1)}(\bB)\]
(see \ref{1.9}). This still leaves some aspects of \cite[Problem 50]{M2}
open: are there any implications between attainment in $\hL_{(0)}$ and
$\hL_{(1)}$ sense? Between $\hL_{(0)}$ and $\hL_{(7)}$ sense?

We also carry out the parallel work for the hereditary density. Let us
introduce the respective definitions. The density $d(X)$ of a topological
space $X$ is defined as the minimal size of a dense subset of $X$. The {\em
topological density\/} $d(\bB)$ of a Boolean algebra $\bB$ is the density of
the space $\Ult(\bB)$ of ultrafilters on $\bB$. The {\em algebraic
density\/} (sometimes also called {\em the $\pi$--weight\/}) of a Boolean
algebra $\bB$ is  
\[\pi(\bB)=\min\{|A|:A\subseteq\bB\setminus\{{\bf 0}\}\mbox{ and } (\forall
b\in\bB\setminus\{{\bf 0}\})(\exists a\in A)(a\leq b)\;\}.\]

\begin{definition}
\label{hdbasic}
For an infinite Boolean algebra $\bB$ we let:
\[\begin{array}{l}
\hd^{(+)}_{(0)}(\bB)=\sup\{d(X)^{(+)}: X\mbox{ is a subspace of }\Ult(\bB)
\;\},\\
\hd^{(+)}_{(5)}(\bB)=\sup\{\kappa^{(+)}:\mbox{ there is a left--separated
sequence of length }\kappa\;\},\\
\hd^{(+)}_{(7)}(\bB)=\sup\{\pi(\bB^*)^{(+)}:\bB^*\mbox{ is a homomorphic
image of }\bB\;\},\\ 
\hd^{(+)}_{(8)}(\bB)=\sup\{d(\bB^*)^{(+)}: \bB^*\mbox{ is a homomorphic
image of }\bB\;\}. 
  \end{array}\]
(Again, the superscripts ``$(+)$'' mean that we have two variants for each
cardinal: with and without ``$+$''.)
\end{definition}
Like before, the cardinals mentioned in \ref{hdbasic} correspond to those
listed in \cite[Theorem 16.1]{M2}, and the variants with ``$+$'' reflect the
attainment properties. The known dependencies here are
\[\begin{array}{l}
\hd^{+}_{(5)}(\bB)=\hd_{(5)}(\bB)\quad\Rightarrow\quad
\hd^{+}_{(7)}(\bB)=\hd_{(7)}(\bB)\quad\Rightarrow\\
\hd^{+}_{(0)}(\bB)=\hd_{(0)}(\bB)\quad\Rightarrow\quad
\hd^{+}_{(8)}(\bB)=\hd_{(8)}(\bB)\end{array}\]
and 
\[\hd_{(0)}(\bB)=\hd_{(0)}^+(\bB)\mbox{ is a regular cardinal }\quad
\Rightarrow\quad \hd_{(5)}=\hd_{(5)}^+(\bB)\]
(and Monk \cite[Problem 54]{M2} asked for a complete description of
dependencies). Like for $\hL$, if $\hd(\bB)$ is a strong limit cardinal or
if it has countable cofinality, then $\hd_{(5)}(\bB)<\hd^+_{(5)}(\bB)$ (see
Juh\'asz \cite[4.2, 4.3]{Ju80}).  

In \ref{233forhd} we note that if $\hd(\bB)$ is a singular cardinal such
that $2^{\cf(\hd(\bB))}<\hd(\bB)$, then $\hd_{(8)}^+(\bB)=\hd_{(7)}^+ 
(\bB)=\hd_{(5)}^+(\bB)=\hd_{(0)}^+(\bB)=(\hd(\bB))^+$. Consequently, GCH
implies that the sups in all equivalent definitions of $\hd$ are attained at
singular cardinals. Then, in section 4, we show that, consistently, there is
a Boolean algebra $\bB$ such that 
\[\hd_{(5)}(\bB)<\hd^+_{(5)}(\bB)\quad\mbox{ and }\quad\hd^+_{(7)}(\bB)=
\hd_{(7)}(\bB)\]
(see \ref{2.6}). This still leaves several aspects of \cite[Problem 54]{M2}
open. 
 
Finally, in the last section of the paper we show that (if we start with the
right cardinals $\mu,\lambda$, $\cf(\lambda)<\lambda$) adding a $\mu$--Cohen
real produces a Boolean algebra $\bB$ such that $\hL_{(7)}^+(\bB)=
\hd_{(5)}^+(\bB)=s^+(\bB)=\lambda$ (put \ref{Prop3.4}, \ref{Thm3.5}
together). This result is of interest as it shows how easily we may have
algebras in which the three cardinal functions do not attain their
supremums. (But of course there is the semi-ZFC result of \cite[Theorem
4.2]{Sh:641}.) 
\bigskip

\noindent{\bf Notation:}\qquad Our notation is standard and compatible with
that of classical textbooks on set theory (like Jech \cite{J}) and Boolean
algebras (like Monk \cite{M1}, \cite{M2}). However in forcing considerations
we keep the older tradition that
\begin{center}
{\em
the stronger condition is the greater one.
}
\end{center}
Let us list some of our notation and conventions.

\begin{enumerate}
\item A name for an object in a forcing extension is denoted with a dot
above (like $\dot{X}$) with one exception: the canonical name for a generic
filter in a forcing notion $\bP$ will be called $\Gamma_\bP$. For a
$\bP$--name $\dot{X}$ and a $\bP$--generic filter $G$ over $\V$, the
interpretation of the name $\dot{X}$ by $G$ is denoted by $\dot{X}^G$. 
\item $i,j,\alpha,\beta,\gamma,\delta,\ldots$ will denote ordinals and 
$\kappa,\mu,\lambda,\theta$ will stand for (always infinite) cardinals.
\item For a set $X$ and a cardinal $\lambda$, $[X]^{\textstyle<\lambda}$
stands for the family of all subsets of $X$ of size less than $\lambda$. If 
$X$ is a set of ordinals then its order type is denoted by $\otp(X)$.
\item Sequences of ordinals will be typically called $\sigma,\rho,\eta,\nu$;
the length of a sequence $\rho$ is $\lh(\rho)$; $\nu\vartriangleleft\eta$
means that the sequence $\nu$ in an initial segment of $\eta$. The set
of all sequences of length $\mu$ with values in $\kappa$ will be denoted by
$\mu^{\textstyle\kappa}$. The lexicographic order on sequences of ordinals
will be called $<_{\lex}$.  
\item In Boolean algebras we use $\vee$ (and $\bigvee$), $\wedge$ (and
$\bigwedge$) and $-$ for the Boolean operations. If $\bB$ is a Boolean
algebra, $x\in\bB$ then $x^0=x$, $x^1=-x$. The Stone space of the algebra
$\bB$ (the space of ultrafilters) is called $\Ult(\bB)$. When working in the
Stone space, we identify the algebra $\bB$ with the field of clopen subsets
of $\Ult(\bB)$. 
\item For a subset $Y$ of an algebra $\bB$, the subalgebra of $\bB$
generated by $Y$ is denoted by $\langle Y\rangle_{\bB}$ and the ideal
generated by $Y$ is called $\id_{\bB}(Y)$.
\end{enumerate}

\bigskip

\noindent{\bf Acknowledgements:}\qquad We would like to thank the referee for
valuable comments and suggestions. 

\section{Golden Oldies: the use of [Sh:233]}
In this section we recall how \cite{Sh:233} applies to the attainment
problems. The proofs of \ref{233} and \ref{233ons} were presented in
\cite{Sh:233}, but we recall them here, as we have an impression that those  
beautiful results went somehow unnoticed. Also, as the results of sections 3
and 4 complement the consequences of \cite[Lemma 5.1]{Sh:233} presented
here, it may be convenient for the reader to have all the proofs presented
as well.  

\begin{hypothesis}
\label{hypoldies}
Let $\mu,\lambda$ be cardinals, and $\bar{\chi}=\langle\chi_i:i<\cf(\lambda)
\rangle$ be an increasing sequence of regular cardinals such that 
\[\cf(\lambda)<\mu=\left(2^{\cf(\lambda)}\right)^+<\lambda=\sup_{i<
\cf(\lambda)}\chi_i \quad\mbox{ and }\quad \mu<\chi_0.\]
\end{hypothesis}

\begin{theorem}
[See {\cite[Lemma 5.1]{Sh:233}}]
\label{233}
Let $X$ be a topological space with a basis $\cB$ consisting of clopen
sets. Suppose that $\Phi$ is a function assigning cardinal numbers to
subsets of $X$ such that $\Phi(X)\geq\lambda$ and 
\begin{enumerate}
\item[(i)] $\Phi(A)\leq\Phi(A\cup B)\leq\Phi(A)+\Phi(B)+\aleph_0$ for
$A,B\subseteq X$, 
\item[(ii)] for each closed set $Y\subseteq X$ such that $\Phi(Y)\geq
\lambda$ and for $i<\cf(\lambda)$, there are $\langle u_\alpha: \alpha<\mu
\rangle\subseteq\cB$ and $\langle y_\alpha:\alpha<\mu\rangle\subseteq Y$
such that
\begin{enumerate}
\item[(a)] $y_\alpha\in u_\alpha\cap Y$,
\item[(b)] $(\forall v\in\cB)(y_\alpha\in v\ \Rightarrow\ \Phi(v\cap Y)\geq
\chi_i)$, 
\item[(c)] $(\forall g:\mu\longrightarrow 2^{\cf(\lambda)})(\exists\alpha,
\beta<\mu)(g(\alpha)=g(\beta)\ \&\ y_\alpha\notin u_\beta)$, 
\end{enumerate}
\item[(iii)] if $\langle A_\alpha:\alpha<\mu\rangle$ is a sequence of
subsets of $X$ such that $\Phi(A_\alpha)\leq\chi_i$ (for $\alpha<\mu$) then
$\Phi(\bigcup\limits_{\alpha<\mu} A_\alpha)\leq\chi_i$.
\end{enumerate}
Then there is a sequence $\langle v_i:i<\cf(\lambda)\rangle\subseteq\cB$
such that 
\[(\forall i<\cf(\lambda))(\Phi(v_i\setminus \bigcup\limits_{j\neq i} v_j)
\geq\chi_i).\] 
\end{theorem}

\begin{proof}
First, by induction on $i<\cf(\lambda)$, we choose families $K_i$ of clopen
subsets of $X$, and sets $D_i\subseteq X$ such that $|K_i|=|D_i|=\mu$. So
suppose that $K_j,D_j$ have been defined for $j<i$. For each $\cU\in
[\bigcup\limits_{j<i} K_j]^{\textstyle{<}\cf(\lambda)}$ such that
$\Phi(X\setminus\bigcup\cU)\geq\lambda$ pick $\langle y_\alpha^\cU:\alpha<
\mu\rangle\subseteq X\setminus\cU$ and $\langle u^\cU_\alpha:\alpha<\mu
\rangle\subseteq\cB$ as guaranteed by (ii) (for $i$ and $Y=X\setminus\bigcup
\cU$). Let $D_i$ consist of all $y^\cU_\alpha$ (for $\cU$ as above and
$\alpha<\mu$); note that $|D_i|=\mu$. Let $K_i$ be a family of clopen sets
such that $|K_i|=\mu$ and for each $\cU$ as above:
\begin{itemize}
\item $u^\cU_\alpha\in K_i$ for all $\alpha<\mu$,
\item if $y^\cU_\alpha\in u^\cU_\alpha\setminus u^\cU_\beta$,
$\alpha,\beta<\mu$, then there is $u\in K_i\cap\cB$ such that $y^\cU_\alpha
\in u\subseteq u^\cU_\alpha\setminus u^\cU_\beta$,
\item if $u\in K_i$ then $X\setminus u\in K_i$. 
\end{itemize}
Let $K=\bigcup\limits_{i<\cf(\lambda)} K_i$ (clearly $|K|=\mu$) and let  
\[Z_i=\{x\in X:\mbox{ if }\{v_\xi:\xi<\cf(\lambda)\}\subseteq K\mbox{ and }
x\in\bigcap\limits_{\xi<\cf(\lambda)} v_\xi\mbox{ then }\Phi(
\bigcap\limits_{\xi<\cf(\lambda)} v_\xi)\geq\chi_i\}.\]

\begin{claim}
\label{cly1}
If $Y\subseteq X$ is a closed set such that $\Phi(Y)\geq\chi_i$, then
$Z_i\cap Y\neq\emptyset$.
\end{claim}

\begin{proof}[Proof of the claim]
Suppose that for each $x\in Y$ we have a sequence $\langle v^x_\xi:\xi<\cf
(\lambda)\rangle\subseteq K$ such that $x\in\bigcap\limits_{\xi<\cf(\lambda
)} v_\xi^x$ and $\Phi(\bigcap\limits_{\xi<\cf(\lambda)}
v_\xi^x)<\chi_i$. There are at most $\mu$ possibilities for such sequences,
so we get a set $W\in [Y]^{\textstyle{\leq}\mu}$ such that 
\[Y\subseteq\bigcup_{x\in W}\bigcap_{\xi<\cf(\lambda)} v^x_\xi.\]
Use the assumption (iii) to conclude that $\Phi(\bigcup\limits_{x\in W}
\bigcap\limits_{\xi<\cf(\lambda)} v^x_\xi)\leq\chi_i$, and next use
(i) to get a contradiction with $\Phi(Y)\geq\lambda$.
\end{proof}

For each $i<\cf(\lambda)$ fix $z_i\in Z_i$.

Now, by induction on $i<\cf(\lambda)$, choose $v_i\in K_i$ and $x_i\in Z_i$
such that 
\begin{enumerate}
\item[$(\alpha)$] $x_i\in v_i\setminus\bigcup\limits_{j<i} v_j$,
$v_i\in\cB$,
\item[$(\beta)$]  $x_j\notin v_i$ for $j<i$,
\item[$(\gamma)$] $z_\vare\notin v_i$ for $i<\vare<\cf(\lambda)$.
\end{enumerate}
Suppose that $x_j,v_j$ have been defined for $j<i$. Let $\cU=\{v_j:j<i\}$
and $Y=X\setminus\bigcup\cU$ (so it is a closed subset of $X$). By
$(\gamma)$, for $\vare>i$ we have $z_\vare\in Y$ and thus $\Phi(Y)\geq
\chi_\vare$ (just look at the definition of $Z_\vare$; remember $X\setminus
v_j\in K$), and hence $\Phi(Y)\geq\lambda$. Consequently, we have sequences
$\langle y^\cU_\alpha:\alpha<\mu\rangle\subseteq D_i$ and $\langle
u^\cU_\alpha:\alpha<\mu\rangle\subseteq K_i$ as chosen before (so they are
as in (ii)). Consider a function $g$ defined on $\mu$ such that 
\[g(\alpha)=u^\cU_\alpha\cap (\{z_\vare:\vare<\cf(\lambda)\}\cup
\{x_j:j<i\}).\]
So by (ii)(c) we find distinct $\alpha,\beta<\mu$ such that $g(\alpha)=
g(\beta)$ and $y^\cU_\alpha\notin u^\cU_\beta$. Then, by the definition of
$K_i$, we find $v_i\in K_i\cap\cB$ such that $y^\cU_\alpha\in v_i\subseteq
u^\cU_\alpha\setminus u^\cU_\beta$. It follows from (ii)(b) that
$\Phi(v_i\cap Y)=\Phi(v_i\setminus\bigcup\limits_{j<i} v_j)\geq\chi_i$. By
claim \ref{cly1} we may pick $x_i\in Z_i\cap v_i\cap Y= Z_i\cap v_i
\setminus\bigcup\limits_{j<i} v_j$. Since, by our choices, $v_i$ is disjoint
from $\{z_\vare:\vare<\cf(\lambda)\}\cup\{x_j:j<i\}$, the inductive step is
complete. 

After the inductive construction is carried out, look at the sequence
$\langle v_i:i<\cf(\lambda)\rangle$. Since $x_i\in Z_i\cap v_i\setminus
\bigcup\limits_{j\neq i} v_j$ we easily conclude that $\Phi(v_i\setminus
\bigcup\limits_{j\neq i} v_j)\geq\chi_i$.
\end{proof}

\begin{corollary}
[See {\cite[3.3., 5.4]{Sh:233}}]
\label{233ons}
If $\bB$ is a Boolean algebra satisfying $s(\bB)=\lambda$ then
$s^+(\bB)=\lambda^+$. 
\end{corollary}

\begin{proof}
Suppose $s(\bB)=\lambda$. Then for each $i<\cf(\lambda)$ we may pick a
discrete set $A_i\subseteq\Ult(\bB)$ of size $\chi_i$. Let
$X=\bigcup\limits_{i<\cf(\lambda)} A_i$ (and the topology of $X$ is the one
inherited from $\Ult(\bB)$) and let $\cB=\{b\cap X: b\in\bB\}$. Finally let
$\Phi(A)=|A|$ for $A\subseteq X$. Note that $X,\cB,\Phi$ clearly satisfy
clauses \ref{233}(i,iii). Suppose that the demand in \ref{233}(ii) fails for 
$i<\cf(\lambda)$ and a closed set $Y\subseteq X$ (so $|Y|=\lambda$). Let 
\[Y^*_i=\{y\in Y:(\forall v\in\cB)(y\in v\ \Rightarrow\ |v\cap Y|\geq
\chi_i)\}.\]

\noindent {\sc Case 1:}\quad $|Y^*_i|<\mu$.\\
Then $|Y\setminus Y^*_i|=\lambda$. For each $y\in Y\setminus Y^*_i$ pick
$v^y\in\cB$ such that $y\in v^y$ and $|v^y\cap Y|<\chi_i$. Consider the
function 
\[F:Y\setminus Y^*_i\longrightarrow{\mathcal P}(Y\setminus Y^*_i):
y\mapsto v^y\cap Y\setminus Y^*_i.\]
By the Hajnal Free Set Theorem (see Hajnal \cite{Ha61}) there is an
$F$--free set $S\subseteq Y\setminus Y^*_i$ of size $\lambda$. Then $y\notin 
F(y')$ for distinct $y,y'\in S$, and thus $v^y\cap S=\{y\}$ for $y\in S$. 
Consequently $S$ is discrete and $s^+(\bB)>\lambda$.  

\noindent {\sc Case 2:}\quad $|Y^*_i|\geq\mu$.\\
For some $j<\cf(\lambda)$ we have $|Y^*_i\cap A_j|\geq\mu$, so we may choose
distinct $y_\alpha\in Y^*_i\cap A_j$ for $\alpha<\mu$. The set $\{y_\alpha:
\alpha<\mu\}$ is discrete (as $A_j$ is so), so we may pick $u_\alpha\in\cB$
such that $(\forall\alpha,\beta<\mu)(y_\alpha\in u_\beta\ \Leftrightarrow\
\alpha=\beta)$. Then $\langle y_\alpha,u_\alpha:\alpha<\mu\rangle$ is as
required in \ref{233}(ii), contradicting our assumption that this clause
fails. 
\medskip

So we may assume that the assumptions of \ref{233} are satisfied, and
therefore we may find $\langle v_i:i<\cf(\lambda)\rangle\subseteq\cB$ such
that $|v_i\setminus\bigcup\limits_{j\neq i}v_j|\geq\chi_i$ for each $i<\cf
(\lambda)$. Then, for every $i<\cf(\lambda)$, there is $\xi(i)<\cf(\lambda)$
such that 
\[|A_{\xi(i)}\cap v_i\setminus\bigcup_{j\neq i}v_j|\geq\chi_i.\]
Let 
\[A=\bigcup_{i<\cf(\lambda)} (A_{\xi(i)}\cap v_i\setminus \bigcup_{j\neq i}
v_j).\]
Clearly $|A|=\lambda$ and easily $A$ is discrete. 
\end{proof}

\begin{theorem}
\label{233forhL}
If $\bB$ is a Boolean algebra satisfying $\hL(\bB)=\lambda$ then
\[\hL_{(0)}^+(\bB)=\hL_{(1)}^+(\bB)=\hL_{(7)}^+(\bB)=\lambda^+.\] 
\end{theorem}

\begin{proof}
So assume $\hL(\bB)=\lambda$. 

If $s^+(\bB)>\lambda$, that is if $\bB$ has an ideal independent sequence of
length $\lambda$, then easily all sups in the equivalent definitions of
$\hL$ are obtained. So we may assume 
\begin{enumerate}
\item[$(\circledast)$] $s^+(\bB)\leq\lambda$ and thus, by \ref{233ons},
$s^+(\bB)<\lambda$. We may also assume that $s^+(\bB)<\chi_0$.
\end{enumerate}
Let $X=\Ult(\bB)$, $\cB=\bB$, and for $Y\subseteq X$ let
\[\Phi(Y)=\sup\{\kappa:\mbox{ there is a right separated sequence in $Y$ of
length }\kappa\;\}.\]
(Recall that in a topological space $Y$, a sequence $\langle y_\xi:\xi<
\kappa\rangle$ is right separated whenever all initial segments of the
sequence are open in the relative topology.) We are going to apply \ref{233}
to $X,\cB,\Phi$, and for that we need to check the assumptions there. 
Clauses (i) and (iii) are obvious, and let us verify \ref{233}(ii). 

Let $i<\cf(\lambda)$ and let $Y\subseteq\Ult(\bB)$ be a closed set such that 
$\Phi(Y)=\lambda$. Let $\langle x_\xi:\xi<\chi^+_i\rangle\subseteq Y$ be a
right separated sequence, and let $b_\xi\in\bB$ be such that $x_\xi\in
b_\xi$ and $x_\zeta\notin b_\xi$ for $\xi<\zeta<\chi_i^+$. Let 
\[Z=\{\xi<\chi_i^+:\cf(\xi)=\chi_i\ \&\ (\exists a\in\bB)(x_\xi\in a\ \&\
\Phi(a\cap Y)<\chi_i)\}.\]

\begin{claim}
\label{cly3}
$Z$ is not stationary in $\chi_i^+$.
\end{claim}

\begin{proof}[Proof of the claim]
Assume $Z$ is stationary. For $\xi\in Z$ pick $a_\xi\in\bB$ such that $x_\xi
\in a_\xi$ and $\Phi(a_\xi\cap Y)<\chi_i$. Note that then for some
$\zeta(\xi)<\xi$ we have 
\[(\forall\vare<\xi) (x_\vare\in a_\xi\ \Rightarrow\ \vare<\zeta(\xi)).\]
By the Fodor lemma, for some $\zeta^*$ the set $Z^*=\{\xi\in Z:\zeta(\xi)=
\zeta^*\}$ is stationary. Now look at the set $Y^*=\{x_\xi:\xi\in Z^*\ \&\
\xi>\zeta^*\}$: we have   
\[(\forall\xi\in Z^*\setminus (\zeta^*+1))((a_\xi\cap b_\xi)\cap Y^*=\{
x_\xi\}).\]
Consequently $Y^*$ is a discrete set of size $\chi_i^+$, contradicting
$(\circledast)$. 
\end{proof}

Thus we may pick an increasing sequence $\langle\xi(\alpha):\alpha<\mu
\rangle$ of ordinals below $\chi_i^+$ such that $\cf(\xi(\alpha))=\chi_i$
and $\xi(\alpha)\notin Z$ (for $\alpha<\mu$). Let $y_\alpha=x_{\xi(\alpha)}$
and $u_\alpha=b_{\xi(\alpha)}$. Then $\langle y_\alpha,u_\alpha:\alpha<\mu
\rangle$ is as required in \ref{233}(ii) (for $Y,i$). 

Consequently we may apply \ref{233} to choose a sequence $\langle v_i:i<
\cf(\lambda)\rangle\subseteq\bB$ such that 
\[(\forall i<\cf(\lambda))(\Phi(v_i\setminus\bigcup_{j\neq i} v_j)\geq
\chi_i).\]
For $i<\cf(\lambda)$ choose a right separated sequence $\langle y^i_\xi:
\xi<\chi_i\rangle\subseteq v_{i+1}\setminus\bigcup\limits_{j\neq i+1}
v_j$. Let $I$ consist of those $b\in\bB$ that for some finite set $W
\subseteq\cf(\lambda)$ and a sequence $\langle\zeta(i):i\in W\rangle\in
\prod\limits_{i\in W}\chi_i$ we have
\[(\forall i<\cf(\lambda))(\forall\xi<\chi_i)(y^i_\xi\in b\ \Rightarrow\
i\in W\ \&\ \xi<\zeta(i)).\]

\begin{claim}
\label{cly2}
$I$ is an ideal in $\bB$ and $\cof(I)=\lambda$. Consequently $\hL^+_{(1)}
(\bB)=\lambda^+$ and hence $\hL^+_{(7)}(\bB)=\lambda^+$.
\end{claim}

\begin{proof}[Proof of the claim]
Plainly, $I$ is an ideal in $\bB$. Suppose that $ A\subseteq I$ is of size
less than $\lambda$, and for $b\in A$ let $W_b\in [\cf(\lambda)]^{\textstyle
{<}\omega}$, $\langle\zeta_b(i):i\in W_b\rangle\in\prod\limits_{i\in W_b}
\chi_i$ witness $b\in I$. Let $i<\cf(\lambda)$ be such that $\chi_i>|A|$ and
let $\sup\{\zeta_b(i):(\exists b\in A)(i\in W_b)\}<\xi<\chi_i$. Take $b\in
\bB$ such that $y^i_\xi\in b$  and $(\forall\zeta<\chi_i)(\xi<\zeta\
\Rightarrow\ y^i_\zeta\notin b)$. Then 
\[y^j_\vare\in b\cap v_{i+1}\quad\Rightarrow\quad j=i\ \&\ \vare\leq\xi,\]
so $v_{i+1}\wedge b\in I$, but it is not included in any member of $Z$.
\end{proof}

Let $Y=\{y^i_\zeta: i<\cf(\lambda)\ \&\ \zeta<\chi_i\}$. 

\begin{claim}
\label{cly4}
$L(Y)=\lambda$, and  consequently $\hL^+_{(0)}(\bB)=\lambda^+$.
\end{claim}

\begin{proof}[Proof of the claim]
For $i<\cf(\lambda)$ and $\xi<\chi_i$, let $U_{i,\xi}$ be an open subset
of $v_{i+1}$ such that 
\[(\forall\zeta<\chi_i)(y^i_\zeta\in U_{i,\xi}\ \Leftrightarrow\
\zeta\leq\xi).\]  
Put $\cU_i=\{U_{i,\xi}:\xi<\chi_i\}$, $\cU=\bigcup\limits_{i<\cf(\lambda
)} \cU_i$. It should be clear that if $\cU'\subseteq \cU_i$ is of size less
than $\chi_i$ then $Y\cap \bigcup\cU'\neq Y\cap \bigcup\cU_i$. Also 
$y^i_\xi\notin \bigcup\cU_j\subseteq v_j$ for $i\neq j$, so we may conclude
that no subfamily of $\cU$ of size less than $\lambda$ covers $Y$, showing
the claim.  
\end{proof}
\end{proof}

\begin{theorem}
\label{233forhd}
If $\hd(\bB)=\lambda$ then $\hd_{(8)}^+(\bB)=\lambda^+$ (and thus also
$\hd_{(0)}^+(\bB)=\hd_{(7)}^+(\bB)=\hd_{(5)}^+(\bB)=\lambda^+$).  
\end{theorem}

\begin{proof}
We may follow like in \ref{233forhL} and use \ref{233} to get our
conclusion. However, an alternative way is to use a result of {\v 
S}apirovski{\u\i} that for every compact space $X$, $\hd(X)\leq s(X)^+$ (see
{\v S}apirovski{\u\i} \cite{Sa74} or Hodel \cite[7.17]{Ho84}). Consequently,
in our situation, $\hd(\bB)=s(\bB)$ and by \ref{233ons} we conclude that
$s^+(\bB)=\lambda^+$. But this implies that there is a homomorphic image
$\bB^*$ of $\bB$ with the cellularity $c(\bB^*)=\lambda$ (see \cite[Theorem
3.25 and p. 175]{M2}). Clearly $d(\bB^*)\geq c(\bB^*)$, so we get our
conclusion. 
\end{proof}

\section{Some combinatorics}
Arguments based on the $\Delta$--lemma are very important in forcing
considerations. The result quoted below is a variant of the $\Delta$--lemma
and in various forms was presented, proved and developed in \cite[\S
6]{Sh:430}, \cite[\S 6]{Sh:513} and \cite[\S 7]{Sh:620}. 

\begin{lemma}
[see {\cite[6.1]{Sh:513}}]
\label{0.A}
Assume that:
\begin{enumerate}
\item[(i)]    $\sigma,\theta$ are regular cardinals and $\kappa$ is a
cardinal, 
\item[(ii)]   $(\forall\alpha<\sigma)(|\alpha|^\kappa<\sigma)$,
\item[(iii)]  $\cD$ is a $\sigma$--complete filter on $\theta$ containing
all co-bounded subsets of $\theta$,
\item[(iv)]   $\langle\beta^\alpha_\varepsilon:\varepsilon<\kappa\rangle$ is 
a sequence of ordinals (for $\alpha<\theta$),
\item[(v)]    $X\subseteq\theta$ is such that $X\neq\emptyset\mod\cD$.
\end{enumerate}
Then there are a sequence $\langle\beta^*_\varepsilon:\varepsilon<\kappa
\rangle$ and a set $w\subseteq\kappa$ such that:
\begin{enumerate}
\item[(a)] $(\forall\varepsilon\in\kappa\setminus w)(\sigma\leq\cf(
\beta^*_\varepsilon)\leq\theta)$, 
\item[(b)] the set 
\[\begin{array}{ll} 
B\stackrel{\rm def}{=}\{\alpha\in X:&\mbox{if }\varepsilon\in w\mbox{ then }
\beta^\alpha_\varepsilon=\beta^*_\varepsilon,\\
&\mbox{if }\varepsilon\in\kappa\setminus w\mbox{ then }\sup\{\beta^*_\zeta:
\zeta<\kappa,\ \beta^*_\zeta<\beta^*_\varepsilon\}<\beta^\alpha_\varepsilon< 
\beta^*_\varepsilon\}
  \end{array}\]
is not\ $\emptyset$\ modulo the filter $\cD$,
\item[(c)] if $\beta_\varepsilon^\prime<\beta^*_\varepsilon$ (for
$\varepsilon\in\kappa\setminus w$) then
\[\{\alpha\in B: (\forall\varepsilon\in\kappa\setminus w)(
\beta^\prime_\varepsilon<\beta^\alpha_\varepsilon)\}\neq\emptyset\mod\cD.\]
\end{enumerate}
\end{lemma}

The above version of the $\Delta$--lemma will have multiple use in our
proofs in the next two sections. Among others, it will be applied to filters
given by \ref{1.2}, \ref{2.1} below.   

\begin{lemma}
\label{1.2}
Suppose that $\bB$ is a Boolean algebra generated by $\langle x_\xi:\xi<\chi
\rangle$. Let $I\subseteq\bB$ be an ideal with $\cof(I)=\lambda$ and let
$\aleph_0<\mu<\lambda$. Then there are a regular cardinal $\theta\in [\mu,
\lambda]$, a $(<\theta)$--complete filter $\cD$ on $\theta$ and a sequence
$\langle a_\alpha:\alpha<\theta\rangle\subseteq I$ such that 
\begin{enumerate}
\item[$(*_1)$] all co-bounded subsets of $\theta$ are in the filter $\cD$,
and for every $b\in I$:
\[\{\alpha<\theta:a_\alpha\leq b\}=\emptyset\mod\cD,\]
\item[$(*_2)$] for each $\alpha<\theta$, $a_\alpha\notin\id_{\bB}(\{a_\beta:
\beta<\alpha\})$, 
\item[$(*_3)$] every $a_\alpha$ (for $\alpha<\theta$) is of the form
\[a_\alpha=\bigwedge\limits_{\ell<n} x^{t(\alpha,\ell)}_{\xi(\alpha,\ell)}
\qquad\quad\mbox{(where $n<\omega$, $\xi(\alpha,\ell)<\chi$, $t(\alpha,\ell)
<2$).}\] 
\end{enumerate}
\end{lemma}

\begin{proof} 
It is basically like \cite[2.2, 2.3]{Sh:479}, but for reader's convenience
we present the proof fully.

\begin{claim}
\label{1.2.1}
Assume $\mu_0<\lambda$. Then there are a regular cardinal $\theta\in [\mu_0, 
\lambda]$ and a set $Y\in [I]^{\textstyle\theta}$, such that 
\[(\forall Z\in [I]^{\textstyle{<}\theta})(\exists b\in Y)(\forall a\in Z)
(b\not\leq a).\] 
\end{claim}

\begin{proof}[Proof of the claim]
Assume not. By induction on $|Y|$ we show that then 
\begin{enumerate}
\item[$(\circledast)$] \quad if $Y\in [I]^{\textstyle{\leq}\lambda}$ then
there is $Y^*\subseteq I$ such that $|Y^*|=\mu_0$ and 
\[(\forall b\in Y)(\exists a\in Y^*)(b\leq a).\] 
\end{enumerate}
If $|Y|\leq\mu_0$, then there is nothing to do.\\
Suppose now that $Y\subseteq I$ and $|Y|>\mu_0$ is a regular cardinal. Then,
using the assumption that the claim fails, we may find a set $Z\subseteq I$
such that $|Z|<|Y|$ and $(\forall b\in Y)(\exists a\in Z)(b\leq a)$. Now
apply the induction hypothesis to $Z$ and get a set $Z^*\subseteq I$ of size
$\mu_0$ such that $(\forall a\in Z)(\exists c\in Z^*)(a\leq c)$ -- clearly
the set $Z^*$ works for $Y$ too.\\
So suppose now that $Y\subseteq I$ and $|Y|$ is a singular cardinal
$>\mu_0$. Let $Y=\bigcup\limits_{\xi<\cf(|Y|)} Y_\xi$, where $|Y_\xi|<|Y|$
(for $\xi<\cf(|Y|)$). For each $\xi$ apply the inductive hypothesis to get 
$Y^*_\xi\subseteq I$ such that $|Y^*_\xi|=\mu_0$ and $(\forall b\in Y_\xi)
(\exists a\in Y^*_\xi)(b\leq a)$. Put $Y^+=\bigcup\limits_{\xi<\cf(|Y|)}
Y^*_\xi$ and note that $|Y^+|\leq\cf(|Y|)\cdot\mu_0<|Y|$. Again, apply the
inductive hypothesis $(\circledast)$, this time to $Y^+$, to get the
respective $Y^*$ and note that it works for $Y$ too. 
\medskip

To finish the proof of the claim note that the statement in $(\circledast)$
contradicts the assumption that $\mu_0<\lambda=\cof(I)$. 
\end{proof}

If a set $Y\subseteq I$ is given by \ref{1.2.1} for $I,\mu_0,\theta$ then we 
say that it is {\em temporarily $(I,\mu_0,\theta)$--good}. 

\begin{claim}
\label{1.2.2}
Suppose that $Y\subseteq I$ is temporarily $(I,\mu,\theta)$--good, $\kappa< 
|Y|$. Assume $Y=\bigcup\limits_{\xi<\kappa} Y_\xi$. Then for some $\xi<
\kappa$ the set $Y_\xi$ is temporarily $(I,\mu,\theta)$--good. 
\end{claim}

\begin{proof}[Proof of the claim]
Suppose that $Y=\bigcup\limits_{\xi<\kappa} Y_\xi$, $\kappa<|Y|$ and no
$Y_\xi$ is temporarily $(I,\mu,\theta)$--good. For $\xi<\kappa$ choose
$Z_\xi\subseteq I$ such that  $|Z_\xi|<|Y|=\theta$ and 
\[(\forall b\in Y_\xi)(\exists a\in Z_\xi)(b\leq a),\]
and put $Z=\bigcup\limits_{\xi<\kappa} Z_\xi$. Then $Z$ contradicts ``$Y$ is
temporarily $(I,\mu,\theta)$--good''. The claim is shown.  
\end{proof}

Now, let $Y\subseteq I$ be a temporarily $(I,\mu,\theta)$--good set, $\theta
=|Y|$, and let $Y=\{b_\alpha:\alpha<\theta\}$ be an enumeration. Each
$b_\alpha$ can be represented as
\[b_\alpha=\bigvee\limits_{j<j_\alpha}\bigwedge\limits_{\ell<n_\alpha}
x^{t(\alpha,j,\ell)}_{\xi(\alpha,j,\ell)}.\] 
By \ref{1.2.2} we find $n^*,j^*$ and $A\in [\theta]^{\textstyle\theta}$ such
that $(\forall\alpha\in A)(j_\alpha=j^*\ \&\ n_\alpha=n^*)$ and the set
$Y^*=\{b_\alpha:\alpha\in A\}$ is temporarily $(I,\mu,\theta)$--good. For
$j<j^*$ and $\alpha\in A$ let $b^j_\alpha=\bigwedge\limits_{\ell<n^*}
x^{t(\alpha,j,\ell)}_{\xi(\alpha,j,\ell)}$ and let $Y^j=\{b^j_\alpha:\alpha
\in A\}$. We claim that for some $j<j^*$ the set $Y^j$ is temporarily
$(I,\mu,\theta)$--good. If not, then we find $Z_j\subseteq I$ (for $j<j^*$)
such that $|Z_j|<\theta$ and $(\forall\alpha\in A)(\exists a\in Z_j)
(b^j_\alpha\leq a)$. Put 
\[Z=\{a_0\vee\ldots\vee a_{j^*-1}: a_0\in Z_0,\ldots,a_{j^*-1}\in
Z_{j^*-1}\}\]
and note that this set contradicts ``$Y^*$ is temporarily $(I,\mu,
\theta)$--good''. 

So let $j_0<j^*$ be such that the set $Y^{**}\stackrel{\rm def}{=}
\{b^{j_0}_\alpha:\alpha\in A\}$ is temporarily $(I,\mu,\theta)$--good and
let $Y^{**}=\{a_\alpha:\alpha<\theta\}$ be an enumeration. 

For $b\in I$ let $F_b=\{\alpha<\theta:a_\alpha\not\leq b\}$ and let $\cD_0$ 
be the $(<\theta)$--complete filter of subsets of $\theta$ generated by
$\{F_b: b\in I\}$.

First note that if $\kappa<\theta$ and $\langle b_\xi:\xi<\kappa\rangle
\subseteq I$ then (by the choice of $Y^{**}$) we may find $\alpha<\theta$
such that $(\forall\xi<\kappa)(a_\alpha\not\leq b_\xi)$. Consequently
$\bigcap\limits_{\xi<\kappa} F_{b_\xi}\neq\emptyset$ and we may conclude
that $\cD_0$ is a proper filter on $\theta$. Since $\alpha\notin
F_{a_\alpha}$, we get that $\cD_0$ extends the filter of co-bounded subsets
of $\theta$.  

\begin{claim}
\label{clx}
The set $A^+\stackrel{\rm def}{=}\{\alpha<\theta:a_\alpha\in\id_{\bB}(
\{a_\beta:\beta<\alpha\})\}$ does not belong to the filter $\cD_0$. 
\end{claim}

\begin{proof}[Proof of the claim]
Assume toward contradiction that $A^+\in\cD_0$. Thus we have a sequence
$\langle b_\xi:\xi<\kappa\rangle\subseteq I$, $\kappa<\theta$, such that
$\bigcap\limits_{\xi<\kappa}F_{b_\xi}\subseteq A^+$. It follows from the
choice of $Y^{**}$ that $Y^{**}\not\subseteq\id_{\bB}(\{b_\xi:\xi<
\kappa\})$. So let $\alpha<\theta$ be the first such that $a_\alpha\notin
\id_{\bB}(\{b_\xi:\xi<\kappa\})$. This implies that $a_\alpha\in
\bigcap\limits_{\xi<\kappa}F_{b_\xi}\subseteq A^+$, and thus $a_\alpha\in
\id_{\bB}(\{a_\beta:\beta<\alpha\})$. By the minimality of $\alpha$ we have
$\id_{\bB}(\{a_\beta:\beta<\alpha\})\subseteq\id_{\bB}(\{b_\xi:\xi<
\kappa\})$, and we get a contradiction. 
\end{proof}

Take the set $A^+$ from \ref{clx} and let $\cD=\{X\setminus A^+:X\in
\cD_0\}$ . Then the filter $\cD$ and $\langle a_\alpha:\alpha\in\theta
\setminus A^+\rangle$ satisfy the demands $(*_1)$--$(*_3)$ (after taking the
increasing enumeration of $\theta\setminus A^+$).  
\end{proof}

\begin{lemma}
[see {\cite[2.2, 2.3]{Sh:479}}]
\label{2.1}
Suppose $\cf(\lambda)<\lambda$, $\mu<\lambda$. Assume that $\bB$ is a
Boolean algebra generated by $\langle x_\xi:\xi<\chi\rangle$ and
$I\subseteq\bB$ is an ideal such that $\pi(\bB/I)=\lambda$. Then there are a 
regular cardinal $\theta\in [\mu,\lambda]$, a $(<\theta)$--complete filter
$\cD$ on $\theta$ and a sequence $\langle
a_\alpha:\alpha<\theta\rangle\subseteq\bB\setminus I$ such that
\begin{enumerate}
\item[$(\otimes_1)$] the filter $\cD$ contains all co-bounded subsets of
$\theta$ and for every $b\in\bB\setminus I$:
\[\{\alpha<\theta:\ b\leq a_\alpha\mod I\}=\emptyset\mod\cD,\]
\item[$(\otimes_2)$] if $\beta<\alpha<\theta$ then $a_\beta\wedge
(-a_\alpha)\notin I$, 
\item[$(\otimes_3)$] every $a_\alpha$ (for $\alpha<\theta$) is of the form
\[a_\alpha=\bigwedge\limits_{\ell<n} x^{t(\alpha,\ell)}_{\xi(\alpha,\ell)}
\qquad\quad\mbox{(where $n<\omega$, $\xi(\alpha,\ell)<\chi$,
$t(\alpha,\ell)<2$).}\]
\end{enumerate}
\end{lemma}

\begin{proof}
It is an easy modification of \cite[2.2, 2.3]{Sh:479} (and the proof is
fully parallel to that of Lemma \ref{1.2} here). 
\end{proof}

One of the ways of describing Boolean algebras is giving a dense set of
ultrafilters (equivalently: homomorphisms from the algebra into 2). This is
useful when we want to force a Boolean algebra by smaller approximations (see
the forcing notions used in \cite{Sh:479}, \cite{RoSh:599}).

\begin{definition}
\label{0.C}
For a set $w$ and a family $F\subseteq 2^{\textstyle w}$ we define

\noindent $\cl(F)=\{g\in 2^{\textstyle w}: (\forall u\in [w]^{\textstyle
<\omega})(\exists f\in F)(f\rest u=g\rest u)\}$,

\noindent $\bB_{(w,F)}$ is the Boolean algebra generated freely by
$\{x_\alpha:\alpha\in w\}$ except that

if $u_0,u_1\in [w]^{\textstyle <\omega}$ and there is no $f\in F$ such
that $f\rest u_0\equiv 0$, $f\rest u_1\equiv 1$

then $\bigwedge\limits_{\alpha\in u_1} x_\alpha\wedge 
\bigwedge\limits_{\alpha\in u_0} (-x_\alpha)=0$.
\end{definition}

\begin{proposition}
[see {\cite[2.6]{Sh:479}}]
\label{0.D}
Let $F\subseteq 2^{\textstyle w}$. Then:
\begin{enumerate}
\item each $f\in F$ extends (uniquely) to a homomorphism from $\bB_{(w,F)}$
to $\{0,1\}$ (i.e.~it preserves the equalities from the definition of
$\bB_{(w,F)}$), 
\item if $\tau(y_0,\ldots,y_\ell)$ is a Boolean term and $\alpha_0,\ldots,
\alpha_\ell\in w$ are distinct then
\[\begin{array}{l}
\bB_{(w,F)}\models\tau(x_{\alpha_0},\ldots,x_{\alpha_\ell})\neq 0\qquad
\qquad\mbox{ if and only if}\\
(\exists f\in F)(\{0,1\}\models\tau(f(\alpha_0),\ldots,f(\alpha_k))=1),
  \end{array}\]
\item if $w\subseteq w^*$, $F^*\subseteq 2^{\textstyle w^*}$ and
\[(\forall f\in F)(\exists g\in F^*)(f\subseteq g)\quad\mbox{ and }\quad
(\forall g\in F^*)(g\rest w\in\cl(F))\]
then $\bB_{(w,F)}$ is a subalgebra of $\bB_{(w^*,F^*)}$.
\end{enumerate}
\end{proposition}

\begin{remark}
\label{0.E}
Let $F\subseteq 2^{\textstyle w}$. We will use the same notation for a
member $f$ of $F$ and the homomorphism from $\bB_{(w,F)}$ determined by
it. Hence, for a Boolean term $\tau$, a finite set $v\subseteq w$ and $f\in 
F$, we may write $f(\tau(x_\alpha:\alpha\in v))$ etc.
\end{remark}

\begin{proposition}
\label{xxx}
Let $\bB$ be a Boolean algebra.
\begin{enumerate}
\item A sequence $\bar{a}=\langle a_\alpha:\alpha<\kappa\rangle$ of elements
of $\bB$ is: 
\begin{itemize}
\item {\em ideal independent} if and only if for each $\alpha<\kappa$ there
is a homomorphism $f_\alpha:\bB\longrightarrow\{0,1\}$ such that  
\[f_\alpha(a_\alpha)=1\ \mbox{ and }\ (\forall\beta<\kappa)(\alpha\neq\beta\
\Rightarrow\ f_\alpha(a_\beta)=0);\] 
\item {\em left--separated} if and only if for each $\alpha<\kappa$ there is a
homomorphism $f_\alpha:\bB\longrightarrow\{0,1\}$ such that 
\[f_\alpha(a_\alpha)=1\ \mbox{ and }\ (\forall\beta<\kappa)(\alpha<\beta\
\Rightarrow\ f_\alpha(a_\beta)=0);\] 
\item {\em right--separated} if and only if for each $\alpha<\kappa$ there
is a homomorphism $f_\alpha:\bB\longrightarrow\{0,1\}$ such that 
\[f_\alpha(a_\alpha)=1\ \mbox{ and }\ (\forall\beta<\alpha)(f_\alpha(
a_\beta) =0).\] 
\end{itemize}
\item If the algebra $\bB$ is generated by a sequence $\langle x_\xi:\xi<
\chi\rangle$, and there is an ideal independent (left--separated,
right--separated, respectively) sequence of elements of $\bB$ of length
$\kappa$, then there is such a sequence with terms of the form 
\[a_\alpha=\bigwedge\limits_{k<k_\alpha} x^{t(\alpha,k)}_{\xi(\alpha,k)}\]
and where $\xi(\alpha,k)<\chi$, $t(\alpha,k)\in\{0,1\}$ and $\xi(\alpha,k)
\neq\xi(\alpha,k')$ whenever $k<k'<k_\alpha$.
\end{enumerate}
\end{proposition}

\section{Forcing for $\hL$}
In this section we show that consistently there is a Boolean algebra $\bB$
of size $\lambda$ in which there is a strictly increasing
$\lambda$--sequence of ideals but every ideal in $\bB$ is generated by less
than $\lambda$ elements. This answers \cite[Problem 43]{M1} (and thus a part
of \cite[Problem 50]{M2}). The problem if the respective example can be
constructed just from cardinal arithmetic assumptions remains open.

\begin{definition}
\label{1.3}
\begin{enumerate}
\item {\em A good parameter} is a tuple $S=(\mu,\lambda,\bar{\chi})$
such that $\mu,\lambda$ are cardinals satisfying 
\[\mu=\mu^{<\mu}<\cf(\lambda)<\lambda\quad\mbox{and}\quad (\forall\alpha<
\cf(\lambda))(\forall\xi<\mu)(\alpha^\xi<\cf(\lambda)),\]
and $\bar{\chi}=\langle\chi_i: i<\cf(\lambda)\rangle$ is a strictly
increasing sequence of regular cardinals such that $\cf(\lambda)<\chi_i<
\lambda$, $(\forall i<\cf(\lambda))(\chi_i^{<\mu}=\chi_i)$ and $\lambda=
\sup\limits_{i<\cf(\lambda)}\chi_i$. 
\item A good parameter $S=(\mu,\lambda,\bar{\chi})$ is {\em a convenient
parameter} if additionally $\cf(\lambda)=\mu^+$.
\end{enumerate}
\end{definition}

\begin{definition}
\label{1.4}
Let $S=(\mu,\lambda,\bar{\chi})$ be a convenient parameter and let the set
\[\xs\stackrel{\rm def}{=}\{(i,\xi): i<\cf(\lambda)\ \&\ \xi<\chi_i\}\]
be equipped with the lexicographic order $\prec_S$ (i.e., $(i,\xi)\prec_S
(i',\xi')$ if and only if either $i<i'$, or $i=i'$ and $\xi<\xi'$).
\begin{enumerate}
\item We define a forcing notion $\qs$ as follows. 

\noindent{\bf A condition} is a tuple $p=\langle w^p,u^p,\langle
f^p_{i,\xi}:(i,\xi)\in u^p\rangle\rangle$ such that 
\begin{enumerate}
\item[(a)] $w^p\in [\cf(\lambda)]^{\textstyle<\!\mu}$, \quad $u^p\in [\xs]^{
\textstyle<\!\mu}$, 
\item[(b)] $(\forall i\in w^p)((i,0)\in u^p)$ and if $(i,\xi)\in u^p$ then
$i\in w^p$,
\item[(c)] for $(i,\xi)\in u^p$,\quad $f^p_{i,\xi}:u^p\longrightarrow 2$ is
a function such that 
\[(j,\zeta)\in u^p\ \&\ (j,\zeta)\prec_S (i,\xi)\quad\Rightarrow\quad
f^p_{i,\xi}(j,\zeta)=0,\]
and $f^p_{i,\xi}(i,\xi)=1$,
\end{enumerate}
\noindent {\bf the order}\quad is given by:\quad $p\leq q$\qquad if and only 
if 
\begin{enumerate}
\item[$(\alpha)$] $w^p\subseteq w^q$,\quad $u^p\subseteq u^q$, \qquad and 
\item[$(\beta)$]  $(\forall (i,\xi)\in u^p)(f^p_{i,\xi}\subseteq f^q_{i,
\xi})$, \qquad and 
\item[$(\gamma)$] for each $(i,\xi)\in u^q$ one of the following occurs:

{\bf either}\ $f^q_{i,\xi}\rest u^p={\bf 0}_{u^p}$,

{\bf or}\ $i\in w^p$ and for some $\zeta,\vare\leq\chi_i$ we have $(i,\zeta) 
\in u^p$ and $f^q_{i,\xi}\rest u^p= (f^p_{i,\zeta})_\vare$, where $(f^p_{i,
\zeta})_\vare:u^p\longrightarrow 2$ is defined by 
\[(f^p_{i,\zeta})_\vare(j,\gamma)=\left\{\begin{array}{ll}
0                   &\mbox{if } j=i,\ \gamma<\vare,\\
f^p_{i,\zeta}(j,\gamma) &\mbox{otherwise,}
\end{array}\right.\]

{\bf or}\  $i\notin w^p$ and $f^q_{i,\xi}\rest u^p= (f^p_{j,\zeta})_\vare$
for some $(j,\zeta)\in u^p$ and $\zeta,\vare\leq\chi_j$, where $(f^p_{j,
\zeta})_\vare$ is defined as above.
\end{enumerate}
\item We say that conditions $p,q\in\qs$ are {\em isomorphic} if the linear 
orders 
\[(u^p,\prec_S\rest u^p)\ \mbox{ and }\ (u^q,\prec_S\rest u^q)\]
are isomorphic, and if $H:u^p\longrightarrow u^q$ is the
$\prec_S$--isomorphism then: 
\begin{enumerate}
\item[$(\alpha)$] $H(i,\xi)=(j,0)$ if and only if $\xi=0$,
\item[$(\beta)$]  $f^p_{i,\xi}=f^q_{H(i,\xi)}\comp H$ (for $(i,\xi)\in
u^p$). 
\end{enumerate}
In this situation we may call $H$ {\em an isomorphism from $p$ to $q$}.
\end{enumerate}
\end{definition}

\begin{remark}
\label{1.5}
\begin{enumerate}
\item Of course, $\prec_S$ is a well ordering of $\xs$ in the order type
$\lambda$. 
\item The forcing notion $\qs$ is a relative of the one used in 
\cite[\S 7]{RoSh:599}.   
\item There are $\mu$ isomorphism types of conditions in $\qs$ (remember
$\mu^{<\mu}=\mu$). A condition $p\in\qs$ is determined by its isomorphism
type and the set $u^p$.  
\end{enumerate}
\end{remark}

\begin{proposition}
\label{1.6}
Let $S=(\mu,\lambda,\bar{\chi})$ be a convenient parameter. Then $\qs$ is a
$(<\mu)$--complete $\mu^+$--cc forcing notion.
\end{proposition}

\begin{proof}
First we should check that $\qs$ is really a partial order and for this we
have to verify the transitivity of $\leq$. So suppose that $p\leq q$ and
$q\leq r$ and let us justify that $p\leq r$. The only perhaps unclear demand
is clause \ref{1.4}(1$\gamma$). Assume that $(i,\xi)\in u^r$ and
$f^r_{i,\xi}\rest u^p\neq {\bf 0}_{u^p}$ and consider two cases.
\medskip

\noindent {\sc Case 1:}\qquad $i\in w^p$.\\
Then $i\in w^q$ and, by the definition of $\leq$ (clause $(\gamma)$), we may
pick $\zeta\leq\vare\leq\chi_i$ such that $(i,\zeta)\in u^q$ and
$f^r_{i,\xi}\rest u^q=(f^q_{i,\zeta})_\vare$. Again by clause $(\gamma)$,
for some $\zeta',\vare'$ we have $(i,\zeta')\in u^p$ and $f^q_{i,\zeta}\rest
u^p= (f^p_{i,\zeta'})_{\vare'}$. Now look at the definition of the operation 
$(\cdot)_\vare$ -- it should be clear that $f^r_{i,\xi}\rest u^p=(f^p_{i,
\zeta'})_{\vare''}$ for some $\vare''$. 
\medskip

\noindent{\sc Case 2:}\qquad $i\notin w^p$.\\
If $i\in w^q$ then for some $\zeta,\vare$ we have $f^r_{i,\xi}\rest
u^q=(f^q_{i,\zeta})_\vare$ and $f^q_{i,\zeta}\rest u^p=(f^p_{j,\zeta'})_{
\vare'}$ for some $j,\zeta',\vare'$. Now, since $i\notin w^p$ we may write
$f^r_{i,\xi}\rest u^p=(f^q_{i,\zeta})_\vare\rest u^p=(f^p_{j,\zeta'})_{
\vare'}$ and we are done. Suppose now that $i\notin w^q$. Then $f^r_{i,\xi} 
\rest u^q=(f^q_{j,\zeta})_\vare$ (for some $j,\zeta,\vare$) and we ask if
$j\in w^p$. If so, then for some $\zeta',\vare'$ we have $f^q_{j,\zeta}
\rest u^p=(f^p_{j,\zeta'})_{\vare'}$ and hence $f^r_{i,\xi}\rest u^p=
(f^p_{j,\zeta'})_{\vare''}$ (for some $\vare''$). If not (i.e., if $j\notin
w^p$) then like before we easily conclude that $f^r_{i,\xi}\rest u^p=
(f^q_{j,\zeta})_\vare\rest u^p=f^q_{j,\zeta}\rest u^p= (f^p_{j',\zeta'})_{
\vare'}$ (for some $j',\zeta',\vare'$). 
\medskip

Thus $\qs$ is a forcing notion. To check that it is $(<\mu)$--complete
suppose that $\gamma<\mu$ and $\langle p_\alpha:\alpha<\gamma\rangle
\subseteq\qs$ is an increasing sequence of conditions. Put $w^p=
\bigcup\limits_{\alpha<\gamma} w^{p_\alpha}$, $u^p=\bigcup\limits_{\alpha<
\gamma}u^{p_\alpha}$ and for $(i,\xi)\in u^p$ let  
\[f^p_{i,\xi}=\bigcup\{f^{p_\alpha}_{i,\xi}: (i,\xi)\in u^{p_\alpha},\
\alpha<\gamma\}.\]  
Plainly, $\langle w^p,u^p,\langle f^p_{i,\xi}:(i,\xi)\in u^p\rangle\rangle
\in\qs$ is an upper bound to $\langle p_\alpha:\alpha<\gamma\rangle$. 

Now assume that $\cA\subseteq\qs$ is of size $\mu^+$. Since $\mu^{<\mu}=\mu$ 
and $\cf(\lambda)=\mu^+$ we may use $\Delta$--lemma and ``standard
cleaning'' and find conditions $p,q\in\cA$ such that 
\begin{enumerate}
\item[(i)]   $p,q$ are isomorphic (and let $H:u^p\longrightarrow u^q$ be the 
isomorphism),
\item[(ii)]  $H\rest (u^p\cap u^q)$ is the identity on $u^p\cap u^q$,
\item[(iii)] $\sup(w^p\cap w^q)<\min (w^p\setminus w^q)\leq\sup (w^p
\setminus w^q)<\min (w^q\setminus w^p)$.
\end{enumerate}
Now we are going to define an upper bound $r$ to $p,q$. To this end we put
$w^r=w^p\cup w^q$, $u^r=u^p\cup u^q$ and for $(i,\xi)\in u^r$ we define
$f^r_{i,\xi}: u^r\longrightarrow 2$ as follows.
\begin{itemize}
\item If $(i,\xi)\in u^p$, $i\in w^p\cap w^q$ then $f^r_{i,\xi}=
f^p_{i,\xi}\cup (f^q_{H(i,\xi)})_\xi$,
\item if $(i,\xi)\in u^q$, $i\in w^p\cap w^q$ then $f^r_{i,\xi}=
(f^p_{H^{-1}(i,\xi)})_\xi\cup f^q_{i,\xi}$,
\item if $(i,\xi)\in u^p$, $i\in w^p\setminus w^q$ then $f^r_{i,\xi}=
f^p_{i,\xi}\cup f^q_{H(i,\xi)}$, 
\item if $(i,\xi)\in u^q$, $i\in w^q\setminus w^p$ then $f^r_{i,\xi}=
{\bf 0}_{u^p}\cup f^q_{i,\xi}$.
\end{itemize}
\noindent It should be clear that in all cases the functions $f^r_{i,\xi}$ 
are well defined and that they satisfy the demand \ref{1.4}(1c). Hence
$r=\langle w^r,u^r,\langle f^r_{i,\xi}: (i,\xi)\in u^r\rangle\rangle\in\qs$
and one easily checks that it is a condition stronger than both $p$ and
$q$. So we may conclude that $\qs$ satisfies the $\mu^+$--chain condition.
\end{proof}

For a condition $p\in\qs$ let $F^p=\{{\bf 0}_{u^p}\}\cup\{(f^p_{i,
\xi})_\zeta:\xi,\zeta\leq\chi_i,\ (i,\xi)\in u^p\}$, where $(f^p_{i,
\xi})_\zeta:u^p\longrightarrow 2$ is defined like in \ref{1.4}(1$\gamma$):  
\[(f^p_{i,\xi})_\zeta(j,\gamma)=\left\{\begin{array}{ll}
0                   &\mbox{if } j=i,\ \gamma<\zeta,\\
f^p_{i,\xi}(j,\gamma) &\mbox{otherwise.}
\end{array}\right.\]
Further, let $\bB_p$ be the Boolean algebra $\bB_{(u^p,F^p)}$ (as defined in 
\ref{0.C}). Note that $p\leq q$ implies that $\bB_p$ is a subalgebra of
$\bB_q$ (remember \ref{0.D}). Let $\dot{\bB}_S^0$ be a $\qs$--name such that
$\forces_{\qs}$`` $\dot{\bB}_S^0=\bigcup\{\bB_p: p\in\Gamma_{\qs}\}$ '' and
for $(i,\xi)\in\xs$ let $\dot{f}_{i,\xi}$ be a $\qs$--name such that 
\[\forces_{\qs}\mbox{`` }\dot{f}_{i,\xi}=\bigcup\{f^p_{i,\xi}:(i,\xi)\in
u^p,\ p\in\Gamma_{\qs}\}\mbox{ ''.}\]  

\begin{proposition}
\label{1.7}
Assume that $S=(\mu,\lambda,\bar{\chi})$ is a convenient parameter. Then in 
$\V^{\qs}$: 
\begin{enumerate}
\item $\dot{f}_{i,\xi}:\xs\longrightarrow 2$ (for $(i,\xi)\in\xs$) is such
that $\dot{f}_{i,\xi}(i,\xi)=1$ and
\[(\forall (j,\zeta)\in\xs)((j,\zeta)\prec_S (i,\xi)\ \Rightarrow\
\dot{f}_{i,\xi}(j,\zeta)=0).\] 
\item $\dot{\bB}_S^0$ is the Boolean algebra $\bB_{(\xs,\dot{F})}$ (see
\ref{0.C}), where   
\[\dot{F}=\{(\dot{f}_{i,\xi})_\zeta:(i,\xi)\in\xs,\ \xi\leq\zeta\leq
\chi_i\}\] 
and $(\dot{f}_{i,\xi})_\zeta:\xs\longrightarrow 2$ is such that
\[(\dot{f}_{i,\xi})_\zeta(j,\gamma)=\left\{\begin{array}{ll}
0                   &\mbox{if } j=i,\ \gamma<\zeta,\\
\dot{f}_{i,\xi}(j,\gamma) &\mbox{otherwise,}
\end{array}\right.\qquad\mbox{(for $(j,\gamma)\in\xs$)}.\]
\item The sequence $\langle x_{i,\xi}: (i,\xi)\in\xs\rangle$ is
right--separated in $\dot{\bB}_S^0$ (when we consider $\xs$ with the well
ordering $\prec_S$).  
\end{enumerate}
\end{proposition}

\begin{proof} 
Should be clear (for the third clause remember that each $\dot{f}_{i,
\xi}$ extends to a homomorphism from $\dot{\bB}_S^0$ to $\{0,1\}$, see
\ref{0.D}; remember \ref{xxx}). 
\end{proof}

\begin{theorem}
\label{1.8}
Assume $S=(\mu,\lambda,\bar{\chi})$ is a convenient parameter. Then 
\[\forces_{\qs}\mbox{`` there is no ideal $I\subseteq\dot{\bB}^0_S$ such
that $\cof(I)=\lambda$ ''.}\]
\end{theorem}

\begin{proof} 
Let $\dot{I}$ be a $\qs$--name for an ideal in $\dot{\bB}^0_S$, $p\in\qs$,
and suppose that $p\forces_{\qs}\cof(\dot{I})=\lambda$.  

Fix $i<\cf(\lambda)$ for a moment. 

\noindent It follows \ref{1.2} that we may choose $p_i$, $\theta_i$, $n_i$,
$\dcD_i$, $\dot{e}_i$ and $\dot{t}_i$ such that
\begin{enumerate}
\item[$(\alpha)$] $p_i\in\qs$ is a condition stronger than $p$, $\theta_i$
is a regular cardinal, $\chi_i^+<\theta_i<\lambda$ and $n_i\in\omega$,
\item[$(\beta)$]  $\dcD_i$ is a $\qs$--name for a $(<\theta_i)$--complete
filter on $\theta_i$ extending the filter of co-bounded subsets of
$\theta_i$, 
\item[$(\gamma)$] $\forces_{\qs}$`` $\dot{e}_i:\theta_i\times n_i
\longrightarrow\xs$ and $\dot{t}_i:\theta_i\times n_i\longrightarrow 2$ ''; 

for $\alpha<\theta_i$ let $\dot{a}^i_\alpha$ be a $\qs$--name for an element
of $\dot{\bB}^0_S$ such that 
\[\forces_{\qs}\mbox{`` }\dot{a}^i_\alpha=\bigwedge\limits_{\ell<n_i} x^{
\dot{t}_i(\alpha,\ell)}_{\dot{e}_i(\alpha,\ell)}\mbox{ ''},\]
\item[$(\delta)$] $p_i\forces_{\qs}$`` $\dot{a}^i_\alpha\in\dot{I}$ ''\quad
for each $\alpha<\theta_i$, 
\item[$(\vare)$]  $p_i\forces_{\qs}$`` if $b\in\dot{I}$ then $\{\alpha<
\theta_i: \dot{a}^i_\alpha\leq b\}=\emptyset\mod\dcD_i$\quad and

\qquad\quad $\dot{a}^i_\alpha\notin\id_{\dot{\bB}^0_S}(\{\dot{a}^i_\beta:
\beta<\alpha\})$ for each $\alpha<\theta_i$ ''.
\end{enumerate}
For each $\alpha<\theta_i$ choose an antichain $\{p^i_{\alpha,\zeta}:\zeta<
\mu\}$ of conditions stronger than $p_i$, maximal above $p_i$, and such that
each $p^i_{\alpha,\zeta}$ decides the values of $\dot{e}_i(\alpha,\cdot)$,
$\dot{t}_i(\alpha,\cdot)$. Let 
\[p^i_{\alpha,\zeta}\forces_{\qs}\mbox{`` }\dot{e}_i(\alpha,\ell)=e^\zeta_i(
\alpha,\ell)\quad\&\quad\dot{t}_i(\alpha,\ell)=t^\zeta_i(\alpha,\ell)\mbox{
''\qquad(for $\ell<n_i$).}\]
Plainly, we may demand that $i\in w^{p^i_{\alpha,\zeta}}$ and $e^\zeta_i(
\alpha,\ell)\in u^{p^i_{\alpha,\zeta}}$ (for $\alpha<\theta_i$, $\zeta<\mu$,
$\ell<n_i$). 

Suppose now that $G\subseteq\qs$ is a generic filter (over $\V$) such that
$p_i\in G$ and work in $\V[G]$ for a while. Since the filter $\dcD_i^G$ is
$(<\theta_i)$--complete we find ordinals $\dot{\gamma}^G_i<\theta_i$ and
$\dot{\zeta}^G_i<\mu$ such that the set 
\[\begin{array}{ll}
\dot{X}^G_i\stackrel{\rm def}{=}\Big\{\beta<\theta_i:&\dot{\gamma}^G_i\leq
\beta\ \mbox{ and }\ p^i_{\dot{\gamma}^G_i,\dot{\zeta}^G_i}, p^i_{\beta,
\dot{\zeta}^G_i}\in G\ \mbox{ and }w^{p^i_{\dot{\gamma}^G_i,
\dot{\zeta}^G_i}}=w^{p^i_{\beta,\dot{\zeta}^G_i}},\\
&\mbox{the conditions }p^i_{\dot{\gamma}^G_i,\dot{\zeta}^G_i}, p^i_{\beta,
\dot{\zeta}^G_i}\mbox{ are isomorphic,\quad and}\\
&\mbox{if }H:u^{p^i_{\dot{\gamma}^G_i,\dot{\zeta}^G_i}}\longrightarrow u^{ 
p^i_{\beta,\dot{\zeta}^G_i}}\mbox{ is the isomorphism then}\\
&(\forall\ell<n_i)(H(e^{\dot{\zeta}^G_i}_i(\dot{\gamma}^G_i,\ell))=e^{\dot{
\zeta}^G_i}_i(\beta,\ell)\ \&\ t^{\dot{\zeta}^G_i}_i(\dot{\gamma}^G_i,\ell)= 
t^{\dot{\zeta}^G_i}_i(\beta,\ell))\\
&\mbox{and if }j\leq i,\ (j,\xi)\in\xs \mbox{ then}\\
&(j,\xi)\in u^{p^i_{\dot{\gamma}^G_i,\dot{\zeta}^G_i}}\ \Leftrightarrow\
(j,\xi)\in u^{p^i_{\beta,\dot{\zeta}^G_i}}\ \Big\} 
  \end{array}\]
is not $\emptyset$ modulo $\dcD^G_i$ (remember that in $\V[G]$ we still have
$\cf(\lambda)^{<\mu}=\cf(\lambda)$ and $\chi_i^{<\mu}=\chi_i$). Let
$\dot{\delta}^G_i=\otp(u^{p^i_{\dot{\gamma}^G_i,\dot{\zeta}^G_i}},\prec_S)$ 
and for $\alpha\in\dot{X}^G_i$ let $\langle s^{\alpha,i}_\vare:\vare<
\dot{\delta}^G_i\rangle$ be the $\prec_S$--increasing enumeration of
$u^{p^i_{\alpha,\dot{\zeta}^G_i}}$. Apply Lemma \ref{0.A} to $\mu^+$,
$\theta_i$, $\dot{\delta}^G_i$, $\dcD^G_i$ and $\langle s^{\alpha,i}_\vare:
\vare<\dot{\delta}^G_i\rangle$ here standing for $\sigma,\theta,\kappa,\cD$
and $\langle \beta^\alpha_\vare:\vare<\kappa\rangle$ (respectively) there. 
(Remember $\prec_S$ is a well ordering of $\xs$ in the order type
$\lambda$.) So we find a sequence $\langle s^{*,i}_\vare:\vare<
\dot{\delta}^G_i\rangle\subseteq\xs$ and a set $\dot{v}^G_i\subseteq
\dot{\delta}^G_i$ such that 
\begin{enumerate}
\item[(i)]   $(\forall\vare\in\dot{\delta}^G_i\setminus\dot{v}^G_i)(\mu^+
\leq\cf(\{s\in\xs:s\prec_S s^{*,i}_\vare\},\prec_S)\leq\theta_i)$, 
\item[(ii)]  the set 
\[\begin{array}{ll}
\dot{B}^G_i\stackrel{\rm def}{=}\Big\{\beta\in\dot{X}^G_i:&\mbox{if }\vare
\in\dot{v}^G_i\mbox{ then } s^{\beta,i}_\vare=s^{*,i}_\vare,\ \ \mbox{
and}\\ 
&\mbox{if }\vare\in\dot{\delta}^G_i\setminus \dot{v}^G_i\mbox{ then}\\
&\sup_{\prec_S}\{s^{*,i}_\zeta:\zeta<\dot{\delta}^G_i,\ s^{*,i}_\zeta
\prec_S s^{*,i}_\vare\}\prec_S s^{\beta,i}_\vare\prec_S s^{*,i}_\vare\ 
\Big\} 
  \end{array}\]
is not $\emptyset$ modulo the filter $\dcD^G_i$,
\item[(iii)] if $s'_\vare\prec_S s^{*,i}_\vare$ for $\vare\in
\dot{\delta}^G_i\setminus \dot{v}^G_i$ then 
\[\big\{\beta\in\dot{B}^G_i:(\forall\vare\in\dot{\delta}^G_i\setminus
\dot{v}^G_i)(s'_\vare\prec_S s^{\beta,i}_\vare)\big\}\neq\emptyset\mod
\dcD^G_i.\]  
\end{enumerate}
As there was no special role assigned to $\dot{\gamma}^G_i$ (other than
determining the order type of a condition) we may assume that
$\dot{\gamma}^G_i\in\dot{B}^G_i$.  

Now we go back to $\V$ and we choose a condition $q_i\in\qs$, ordinals
$\gamma_i,\zeta_i,\delta_i$, a set $v_i$ and a sequence $\langle
s^{*,i}_\vare:\vare<\delta_i\rangle\subseteq\xs$ such that $q_i\geq 
p^i_{\gamma_i,\zeta_i}$ and $q_i$ forces that these objects have the 
properties listed in (i)--(iii) above. Note that if some condition stronger
than $q_i$ forces that $\beta\in\dot{B}_i$, then any condition stronger than
both $q_i$ and $p^i_{\beta,\zeta_i}$ does so. Then the conditions
$p^i_{\beta,\zeta_i}$ and $p^i_{\gamma_i,\zeta_i}$ are isomorphic and the
isomorphism is the identity on $u^{p^i_{\beta,\zeta_i}}\cap u^{p^i_{
\gamma_i,\zeta_i}}$, and it preserves $e^{\zeta_i}_i,t^{\zeta_i}_i$. Also
then $w^{p^i_{\beta,\zeta_i}}=w^{p^i_{\gamma_i,\zeta_i}}$ and $u^{p^i_{
\beta,\zeta_i}}\cap (\{j\}\times\chi_j)=u^{p^i_{\gamma_i,\zeta_i}}\cap
(\{j\}\times\chi_j)$ for $j\leq i$. In this situation we will use $\langle
s^{\beta,i}_\vare:\vare<\delta_i\rangle$ to denote the $\prec_S$--increasing
enumeration of $u^{p^i_{\beta,\zeta_i}}$ (and so $s^{\beta,i}_\vare=
s^{*,i}_\vare$ for $\vare\in v_i$, and $\sup_{\prec_S}\{s^{*,i}_\zeta:\zeta<
\delta_i,\ s^{*,i}_\zeta\prec_S s^{*,i}_\vare\}\prec_S s^{\beta,i}_\vare
\prec_S s^{*,i}_\vare$ for $\vare\in\delta\setminus v_i$).

\begin{claim}
\label{1.8.1}
If $j\leq i<\cf(\lambda)$, $\ell<n_i$ and $e^{\zeta_i}_i(\gamma_i,\ell)=
(j,\vare)$ (for some $\vare$) then $t^{\zeta_i}_i(\gamma_i,\ell)=1$. 
\end{claim}

\begin{proof}[Proof of the claim]
Suppose that the claim fails for some $j_0\leq i$, $\vare_0<\chi_{j_0}$ and
$\ell_0<n_i$ (i.e., $t^{\zeta_i}_i(\gamma_i,\ell_0)=0$ and $e^{\zeta_i}_i(
\gamma_i,\ell_0)=(j_0,\vare_0)$). Choose $\alpha$ such that
$\gamma_i<\alpha<\theta_i$ and letting $r_1=p^i_{\gamma_i,\zeta_i}$,
$r_2=p^i_{\alpha,\zeta_i}$ we have: 
\begin{itemize}
\item the conditions $r_1,r_2$ are isomorphic and if $H$ is the isomorphism 
from $r_1$ to $r_2$ then $H(e^{\zeta_i}_i(\gamma_i,\ell))=e^{\zeta_i}_i(
\alpha,\ell)$ and $t^{\zeta_i}_i(\gamma_i,\ell)=t^{\zeta_i}_i(\alpha,\ell)$
(for $\ell<n_i$), 
\item $w^{r_1}=w^{r_2}$ and the isomorphism $H$ is the identity on $u^{r_1}
\cap u^{r_2}$,
\item $(j,\xi)\prec_S H(j,\xi)$ for $(j,\xi)\in u^{r_1}\setminus u^{r_2}$,
and 
\item if $j\leq i$, $(j,\xi)\in\xs$ then $(j,\xi)\in u^{r_1}\
\Leftrightarrow\ (j,\xi)\in u^{r_2}$.  
\end{itemize}
Why is the choice possible? Let $G\subseteq\qs$ be generic over $\V$ such
that $q_i\in G$. It follows from clauses (ii), (iii) that we may find
$\alpha\in\dot{B}^G_i\setminus (\gamma_i+1)$ such that $(\forall\vare\in
\delta_i\setminus v_i)(s^{\gamma_i,i}_\vare\prec_S s^{\alpha,i}_\vare)$. 
Then the two ordinals $\gamma_i,\alpha$ have the required properties in
$\V[G]$, and hence clearly in $\V$ too. 

Next we let $w^r=w^{r_1}=w^{r_2}$, $u^r=u^{r_1}\cup u^{r_2}$ and for
$(j,\xi)\in u^r$ we define $f^r_{j,\xi}:u^r\longrightarrow 2$ as follows.
\begin{itemize}
\item If $(j,\xi)\in u^{r_1}\cap u^{r_2}$ then $f^r_{j,\xi}= f^{r_1}_{
j,\xi}\cup f^{r_2}_{j,\xi}$,
\item if $(j,\xi)\in u^{r_1}\setminus u^{r_2}$ then $f^r_{j,\xi}=f^{
r_1}_{j,\xi}\cup f^{r_2}_{H(j,\xi)}$,
\item if $(j,\xi)\in u^{r_2}\setminus u^{r_1}$ then $f^r_{j,\xi}=(f^{
r_1}_{H^{-1}(j,\xi)})_\xi\cup f^{r_2}_{j,\xi}$.
\end{itemize}
Check that the functions $f^r_{j,\xi}$ are well defined and that 
\[r=\langle w^r,u^r,\langle f^r_{j,\xi}:(j,\xi)\in u^r\rangle\rangle\in\qs\]
is a condition stronger than $r_1,r_2$. Let $\tau_1=
\bigwedge\limits_{\ell<n_i}x^{t^{\zeta_i}_i(\gamma_i,\ell)}_{e^{\zeta_i}_i(
\gamma_i,\ell)}$ and $\tau_2=\bigwedge\limits_{\ell<n_i}x^{t^{\zeta_i}_i(
\alpha,\ell)}_{e^{\zeta_i}_i(\alpha,\ell)}$. Suppose that $(j,\xi)\in u^r$
and $\xi\leq\zeta<\chi_j$. If $j\leq i$ then $(\{j\}\times\chi_j)\cap
u^{r_1}=(\{j\}\times\chi_j)\cap u^{r_2}$ and therefore $(f^r_{j,\xi})_\zeta(
\tau_1)=(f^r_{j,\xi})_\zeta(\tau_2)$. If $j>i$ then necessarily $(f^r_{j,
\xi})_\zeta(j_0,\xi_0)=0$, so $(f^r_{j,\xi})_\zeta(\tau_1)=(f^r_{j,
\xi})_\zeta(\tau_2)=0$. Consequently $\bB_r\models\tau_1=\tau_2$ and hence
$r\forces\dot{a}^i_{\gamma_i}=\dot{a}^i_{\alpha}$, contradicting clause
$(\vare)$ (and so finishing the proof of the claim). 
\end{proof}

Take $n<\omega$, $\delta<\mu$, $v\subseteq\delta$ and an unbounded set $Y
\subseteq\cf(\lambda)$ such that for $i,j\in Y$:
\begin{itemize}
\item $n_i=n$, $\delta_i=\delta$, $v_i=v$, and
\item the conditions $p^i_{\gamma_i,\zeta_i},p^j_{\gamma_j,\zeta_j}$ are
isomorphic, and the isomorphism maps $e^{\zeta_i}_i(\gamma_i,\cdot)$ and
$t^{\zeta_i}_i(\gamma_i,\cdot)$ onto $e^{\zeta_j}_j(\gamma_j,\cdot),
t^{\zeta_j}_j(\gamma_j,\cdot)$, respectively.
\end{itemize}
Now apply Lemma \ref{0.A} to find a sequence $\langle s_{*,\vare}:\vare<
\delta\rangle\subseteq\xs\cup\{(\cf(\lambda),0)\}$ and a set $v^*\subseteq
\delta$ such that 
\begin{enumerate}
\item[(a)] $(\forall\vare\in\delta\setminus v^*)(\cf(\{s\in\xs:s\prec_S
s_{*,\vare}\},\prec_S)=\mu^+)$,
\item[(b)] the set 
\[\begin{array}{ll}
C\stackrel{\rm def}{=}\Big\{i\in Y:&\mbox{if }\vare\in v^*\mbox{ then }
s^{*,i}_\vare=s_{*,\vare},\ \ \mbox{ and}\\ 
&\mbox{if }\vare\in\delta\setminus v^*\mbox{ then}\\
&\sup_{\prec_S}\{s_{*,\zeta}:\zeta<\delta,\ s_{*,\zeta}\prec_S s_{*,\vare} 
\}\prec_S s^{*,i}_\vare\prec_S s_{*,\vare}\Big\}
  \end{array}\]
is unbounded in $\cf(\lambda)$,
\item[(c)] if $s'_\vare\prec_S s_{*,\vare}$ for $\vare\in\delta\setminus
v^*$, then the set  
\[\{i\in C:(\forall\vare\in\delta\setminus v^*)(s'_\vare\prec_S
s^{*,i}_\vare)\}\] 
is unbounded in $\cf(\lambda)$.
\end{enumerate}
[So $\sigma,\theta,\kappa,\cD$ and $\langle\langle\beta^\alpha_\vare:\vare<
\kappa\rangle:\alpha<\theta\rangle$  in \ref{0.A} correspond to
$\cf(\lambda)=\mu^+$, $\delta^*$ and the filter of co-bounded subsets of
$\cf(\lambda)$ and $\langle\langle s^{*,i}_\vare:\vare<\delta\rangle: i<
\cf(\lambda)\rangle$ here.]

Next we use clauses (c), (a) and (iii), (i) to choose inductively a set
$C^+\subseteq C$ of size $\cf(\lambda)$ and ordinals $\alpha_i<\theta_i$
(for $i\in C^+$) such that for every $i\in C^+$:
\begin{enumerate}
\item[(d)] if $\vare\in\delta\setminus v^*$ then for all $j\in C^+\cap i$
and $\zeta<\delta$ we have
\[s^{*,j}_\zeta\prec_S s_{*,\vare}\ \Rightarrow\ s^{*,j}_\zeta\prec_S
s^{*,i}_\vare\quad\mbox{ and }\quad s^{\alpha_j,j}_\zeta\prec_S s_{*,\vare}\ 
\Rightarrow\ s^{\alpha_j,j}_\zeta\prec_S s^{*,i}_\vare,\]
\item[(e)] some condition stronger than $q_i$ forces that $\alpha_i\in
\dot{B}_i$ (see clause (ii) earlier), 
\item[(f)] if $\vare\in\delta\setminus v$ then for all $j\in C^+\cap i$ and
$\zeta<\delta$ we have 
\[s^{*,j}_\zeta\prec_S s^{*,i}_\vare\ \Rightarrow\ s^{*,j}_\zeta\prec_S
s^{\alpha_i,i}_\vare\quad\mbox{ and }\quad s^{\alpha_j,j}_\zeta\prec_S
s^{*,i}_\vare\ \Rightarrow\ s^{\alpha_j,j}_\zeta\prec_S
s^{\alpha_i,i}_\vare,\]
\item[(g)] if $\vare\in v^*$, $s_{*,\vare}=(j,\zeta)$ then $j<\min(C^+)$.
\end{enumerate}
Note that then 
\[i,j\in C^+\ \&\ \zeta,\vare<\delta\ \&\ s^{\alpha_j,j}_\zeta=s^{\alpha_i,
i}_\vare\quad\Rightarrow\quad \vare=\zeta\in v\cap v^*.\]
So $\langle\langle s^{\alpha_i,i}_\vare:\vare<\delta\rangle:i\in C^+\rangle$
is a $\Delta$--system of sequences with the heart $\langle s_{*,\vare}:\vare
\in v\cap v^*\rangle$. Let $u^*=\{s_{*,\vare}:\vare\in v\cap v^*\}$ and
$w^*=\{j<\cf(\lambda): (j,0)\in u^*\}$.

Pick $i^*\in C^+$ such that $|C^+\cap i^*|=\mu$.

\begin{claim}
\label{1.8.2}
\[q_{i^*}\forces_{\qs}\mbox{`` }(\forall\alpha\in\dot{B}_{i^*})(\exists j_1, 
j_2\in C^+)(\dot{a}^{i^*}_\alpha\leq\dot{a}^{j_1}_{\alpha_{j_1}}\vee
\dot{a}^{j_2}_{\alpha_{j_2}}\ \&\ p_{j_1},p_{j_2}\in\Gamma_{\qs})\mbox{
''}.\] 
\end{claim}

\begin{proof}[Proof of the claim]
We are going to show that for every condition $q\geq q_{i^*}$ and an ordinal
$\alpha<\theta_{i^*}$ such that $q\forces\alpha\in\dot{B}_{i^*}$, there are
a condition $r\geq q$ and ordinals $j_1,j_2\in C^+$ such that   
\[r\forces\mbox{`` }\dot{a}^{i^*}_\alpha\leq\dot{a}^{j_1}_{\alpha_{j_1}}
\vee\dot{a}^{j_2}_{\alpha_{j_2}}\ \&\ p_{j_1},p_{j_2}\in\Gamma_{\qs}\mbox{  
''.}\] 
So suppose $q\geq q_{i^*}$ and $q\forces\alpha\in\dot{B}_{i^*}$. We may
assume that $p^{i^*}_{\alpha,\zeta_{i^*}}\leq q$ (see the definition of
$\dot{X}_{i^*},\dot{B}_{i^*}$). Choose $j_1\in C^+\cap i^*$ and $j_2\in
C^+\setminus (i^*+1)$ such that 
\[u^q\cap u^{p^{j_1}_{\alpha_{j_1},\zeta_{j_1}}}=u^q\cap
u^{p^{j_2}_{\alpha_{j_2},\zeta_{j_2}}}=u^*\quad\mbox{ and }\quad\sup(w^q)<
\min(w^{p^{j_2}_{\alpha_{j_2},\zeta_{j_2}}}\setminus w^*).\]
(Remember that $\{u^{p^j_{\alpha_j,\zeta_j}}:j\in C^+\}$ forms a
$\Delta$--system with heart $u^*$ and hence $\{w^{p^j_{\alpha_j,\zeta_j}}:
j\in C^+\}$ forms a $\Delta$--system with heart $w^*$.)

To make the notation somewhat simpler let $q^0=p^{i^*}_{\alpha,
\zeta_{i^*}}$, $q^1=p^{j_1}_{\alpha_{j_1},\zeta_{j_1}}$ and
$q^2=p^{j_2}_{\alpha_{j_2},\zeta_{j_2}}$. Note that the conditions $q^0,q^1,
q^2$ are pairwise isomorphic, and the isomorphisms are the identity on the
$u^*$ (which is the common part of any two $u^{q^k}$'s). Put 
\[\tau_0=\bigwedge_{\ell<n}x^{t^{\zeta_{i^*}}_{i^*}(\alpha,\ell)}_{e^{ 
\zeta_{i^*}}_{i^*}(\alpha,\ell)}\quad\mbox{ and }\quad\tau_k=\bigwedge_{
\ell<n}x^{t^{\zeta_{j_k}}_{j_k}(\alpha_{j_k},\ell)}_{e^{\zeta_{j_k}}_{j_k}
(\alpha_{j_k},\ell)}\qquad \mbox{ (for $k=1,2$).}\] 
Thus $\tau_k$ is an element of the algebra $\bB_{q^k}$. Clearly, for
$k,k'<3$, the isomorphism $H^{k,k'}$ from $q^k$ to $q^{k'}$ carries $\tau_k$
to $\tau_{k'}$. 

Now we are going to define a condition $r\in\qs$ stronger than $q,q^1$ and
$q^2$. For this we put $w^r=w^q\cup w^{q^1}\cup w^{q^2}$, $u^r=u^q\cup
u^{q^1}\cup u^{q^2}$ and we define functions $f^r_{i,\xi}:u^r\longrightarrow 
2$ considering several cases. 
\begin{enumerate}
\item If $(i,\xi)\in u^{q^1}$ and $i\in w^*$ then we put $f^r_{i,\xi}=
f^q_{H^{1,0}(i,\xi)}\cup f^{q^1}_{i,\xi}\cup f^{q^2}_{H^{1,2}(i,\xi)}$ (note 
that this includes the case $(i,\xi)\in u^*$). 
\item If $(i,\xi)\in u^{q^1}$, $i\notin w^*$ then we put $f^r_{i,\xi}={\bf
0}_{u^q}\cup f^{q^1}_{i,\xi}\cup {\bf 0}_{u^{q^2}}$.
\item If $(i,\xi)\in u^q\setminus u^*$ then we look at $f^q_{i,\xi}\rest
u^{q^0}$. If it is ${\bf 0}_{u^{q^0}}$ then we let $f^r_{i,\xi}=f^q_{i,\xi}
\cup {\bf 0}_{u^{q^1}}\cup {\bf 0}_{u^{q^2}}$. Otherwise we find $(j,\zeta)
\in u^{q^0}$ and $\zeta\leq\vare\leq\chi_j$ such that $f^q_{i,\xi}\rest
u^{q^0}=(f^{q^0}_{j,\zeta})_\vare$ and if $i\in w^{q^0}$ then $i=j$, and we
define: 
\begin{enumerate}
\item[$(\alpha)$] if $j\in w^*$, $j<i\leq\sup(w^*)$ then $f^r_{i,\xi}=f^q_{i,
\xi}\cup (f^{q^1}_{H^{0,1}(j,\zeta)})_{\chi_j}\cup (f^{q^2}_{H^{0,2}(j,
\zeta)})_{\chi_j}$,
\item[$(\beta)$] if $i=j\in w^*$ then $f^r_{i,\xi}=f^q_{i,\xi}\cup
(f^{q^1}_{H^{0,1}(j,\zeta)})_{\vare^*}\cup(f^{q^2}_{H^{0,2}(j,\zeta)})_{
\vare^*}$, where $\vare^*=\max\{\vare,\xi\}$, 
\item[$(\gamma)$] if $j\in w^*$, $i<j$ then $f^r_{i,\xi}=f^q_{i,\xi}\cup
(f^{q^1}_{H^{0,1}(j,\zeta)})_{\vare}\cup(f^{q^2}_{H^{0,2}(j,\zeta)})_\vare$, 
\item[$(\delta)$] if either $i>\sup(w^*)$ or $j\notin w^*$ then we first
choose $j'\in w^{q^2}$ and $\zeta'\leq\vare'\leq\chi_{j'}$ such that $(j',
\zeta')\in u^{q^2}$ and $(f^{q^2}_{j',\xi'})_{\vare'}(j'',\xi'')=0$ whenever
$(j'',\xi'')\in u^{q^2}$, $j''\in w^*$, and $(f^{q^2}_{j',\xi'})_{\vare'}(
\tau_2)=1$ if possible (under our conditions); next we let $f^r_{i,\xi}=
f^q_{i,\xi}\cup {\bf 0}_{u^{q^1}}\cup(f^{q^2}_{j',\zeta'})_{\vare'}$. 
\end{enumerate}
\item If $(i,\xi)\in u^{q^2}\setminus u^*$, $i\in w^*$ then we let
$f^r_{i,\xi}=(f^q_{H^{2,0}(i,\xi)})_\xi\cup (f^{q^1}_{H^{2,1}(i,\xi)})_\xi
\cup f^{q^2}_{i,\xi}$.  
\item If $(i,\xi)\in u^{q^2}$, $i\notin w^*$ then we put $f^r_{i,\xi}={\bf
0}_{u^q}\cup {\bf 0}_{u^{q^1}}\cup f^{q^2}_{i,\xi}$.  
\end{enumerate}
It should be a routine to check that in all cases the function $f^r_{i,\xi}$ 
is well defined and that $r=\langle w^r,u^r,\langle f^r_{i,\xi}:(i,\xi)\in
u^r\rangle\rangle\in\qs$ is a condition stronger than $q,q^1,q^2$ (and thus
stronger than $p_{j_1}, p_{j_2}$). [Remember that $w^*\subseteq \min(C^+)$,
so for $j\in w^*$ we have $(j,\xi)\in u^{q^0}\ \Leftrightarrow\ (j,\xi)\in
u^{p^{i^*}_{\alpha_{i^*},\zeta_{i^*}}}$ and hence, when checking clause
\ref{1.4}(1c) in Case 1, we may use clauses {\rm (d), (f)} of the choice of
the set $C^+$. They imply that if $(i,\xi)\in u^{q^1}$, $i\in w^*$ then
$(i,\xi)\preceq_S H^{1,0}(i,\xi)\preceq_S H^{1,2}(i,\xi)$. Considering Case
3($\delta$) with $j\notin w^*$, use the fact that $\min(w^{q^0}\setminus
w^*)\geq\sup(w^*)$ (it follows from our choices). Similarly in Case 2
remember $\min(w^{q^1}\setminus w^*)\geq\sup (w^*).$]

We claim that $\bB_r\models \tau_0\leq\tau_1\vee\tau_2$ and for this we have
to show that there is no function $f\in F^r$ with $f(\tau_0)=1$ and
$f(\tau_1)=f(\tau_2)=0$ (see \ref{0.D}). So suppose toward contradiction 
that $f\in F^r$ is such a function. Note that $f$ cannot be ${\bf 0}_{u^r}$
as then the values given to all the terms would be the same (remember they
are isomorphic). So for some $(i,\xi)\in u^r$ and $\xi\leq\vare\leq\chi_i$ 
we have $f=(f^r_{i,\xi})_\vare$. Let us look at all the cases appearing in
the definition of the functions $f^r_{j,\zeta}$'s (we keep labeling as there 
so we do not repeat the descriptions of the cases).
\medskip

\noindent{\sc Case 1:}\quad Clearly $f^r_{i,\xi}(\tau_0)=f^r_{i,\xi}
(\tau_1)$. It follows from the demands {\rm (d), (f)} of the choice of $C^+$ 
that if $i\in w^*$, $(i,\zeta)\in u^{q^0}$, $(i',\zeta')=H^{0,1}(i,\zeta)$, 
then $i'=i$ and $\zeta'\leq\zeta$. Consequently, we may use \ref{1.8.1} to
conclude that  $(f^r_{i,\xi})_\vare(\tau_0)\leq (f^r_{i,\xi})_\vare(
\tau_1)$, what contradicts the choice of $f$.   
\medskip

\noindent{\sc Case 2:}\quad Plainly $(f^r_{i,\xi})_\vare(\tau_0)=(f^r_{i,
\xi})_\vare(\tau_2)$.  
\medskip

\noindent{\sc Case 3$\alpha$:}\quad Note that $f^r_{i,\xi}(\tau_0)=f^r_{i, 
\xi}(\tau_1)$ and, as $j<i\leq\sup(w^*)$, necessarily $i\notin w^{q^0}\cup
w^{q^1}$. Hence easily $(f^r_{i,\xi})_\vare(\tau_0)=(f^r_{i,\xi})_\vare
(\tau_1)$.  
\medskip

\noindent{\sc Cases 3$\beta,\gamma$, 4:}\quad Like in cases 1, 3$\alpha$ we
conclude $(f^r_{i,\xi})_\vare(\tau_0)\leq (f^r_{i,\xi})_\vare(\tau_1)$.  
\medskip

\noindent{\sc Case 3$\delta$:}\quad It follows from the choice of $\zeta',
\vare',j'$ there that $f^r_{i,\xi}(\tau_0)\leq f^r_{i,\xi}(\tau_2)$. If
$i\notin w^{q^0}$ then (as also $i\notin w^{q_2}$) we have 
$f(\tau_0)=f^r_{i,\xi}(\tau_0)$ and $f(\tau_2)=f^r_{i,\xi}(\tau_2)$, so we
are done. If $i\in w^{q^0}$ then $i=j$ and we easily finish by the choice of
$\zeta',\vare',j'$. 
\medskip

\noindent{\sc Case 5:}\quad Clearly $(f^r_{i,\xi})_\vare(\tau_0)=(f^r_{i,
\xi})_\vare(\tau_1)$, a contradiction. 
\medskip

Thus we may conclude that $r\forces$`` $\dot{a}^{i^*}_\alpha\leq
\dot{a}^{j_1}_{\alpha_{j_1}}\vee\dot{a}^{j_2}_{\alpha_{j_2}}$ '', finishing
the proof of the claim.  
\end{proof}

Now we may easily finish the theorem: take a generic filter $G\subseteq\qs$
over $\V$ such that $q_{i^*}\in G$ and work in $\V[G]$. Since the filter
$\dcD^G_{i^*}$ is $(<\theta_{i^*})$--complete and $\cf(\lambda)<\theta_{ 
i^*}$, we find $j_1,j_2\in C^+$ such that $p_{j_1},p_{j_2}\in G$ and 
\[\{\alpha\in\dot{B}^G_{i^*}:(\dot{a}^{i^*}_\alpha)^G\leq (\dot{a}^{j_1}_{
\alpha_{j_1}})^G\vee (\dot{a}^{j_2}_{\alpha_{j_2}})^G\}\neq\emptyset\mod
\dcD^G_{i^*}\] 
(remember $\dot{B}^G_{i^*}\neq\emptyset\mod\dcD^G_{i^*}$ by {\rm (ii)}). But 
then also $(\dot{a}^{j_1}_{\alpha_{j_1}})^G\vee (\dot{a}^{j_2}_{\alpha_{
j_2}})^G\in\dot{I}^G$, so we get a contradiction to clause $(\vare)$.
\end{proof}

\begin{conclusion}
\label{1.9}
It is consistent that there is a Boolean algebra $\bB$ of size $\lambda$
such that there is a right--separated sequence of length $\lambda$ in $\bB$,  
(so $\hL_{(7)}^+(\bB)=\lambda^+$), but there is no ideal $I\subseteq\bB$
with the generating number $\lambda$ (and thus $\hL_{(1)}^+(\bB)=\hL_{(1)}
(\bB)=\lambda$).  
\end{conclusion}

\begin{problem}
Does there exist a Boolean algebra $\bB$ as in \ref{1.9} in semi-ZFC? I.e.,
can one construct such an algebra for $\lambda$ from cardinal arithmetic
assumptions?   
\end{problem}

\section{Forcing for $\hd$}
Here we deal with a problem parallel to the one from the previous section
and related to the attainment question for $\hd$. We introduce a forcing
notion $\Ps$ complementary to $\qs$ and we use it to show that,
consistently, there is a Boolean algebra $\bB$ of size $\lambda$ in which
there is a strictly decreasing $\lambda$--sequence of ideals but every
homomorphic image of $\bB$ has algebraic density less than $\lambda$. This 
gives a partial answer to \cite[Problem 54]{M2}. Again, we do not know if an
example like that can be constructed from cardinal arithmetic assumptions. 

\begin{definition}
\label{2.2}
Let $S=(\mu,\lambda,\bar{\chi})$ be a good parameter (see \ref{1.3}) and let 
$\xs,\prec_S$ be as defined in \ref{1.4}.   
\begin{enumerate}
\item We define a forcing notion $\Ps$ as follows. 

\noindent{\bf A condition} is a tuple $p=\langle w^p,u^p,\langle
f^p_{i,\xi}:(i,\xi)\in u^p\rangle\rangle$ such that 
\begin{enumerate}
\item[(a)] $w^p\in [\cf(\lambda)]^{\textstyle<\!\mu}$, \quad $u^p\in 
[\xs]^{\textstyle<\!\mu}$, 
\item[(b)] $(\forall i\in w^p)((i,0)\in u^p)$ and if $(i,\xi)\in u^p$ then
$i\in w^p$,
\item[(c)] for $(i,\xi)\in u^p$,\quad $f^p_{i,\xi}:u^p\longrightarrow 2$ is
a function such that 
\[(j,\zeta)\in u^p\ \&\ (i,\xi)\prec_S (j,\zeta)\quad\Rightarrow\quad
f^p_{i,\xi}(j,\zeta)=0,\]
and $f^p_{i,\xi}(i,\xi)=1$,
\end{enumerate}
\noindent {\bf the order}\quad is given by:\quad $p\leq q$\qquad if and only
if 
\begin{enumerate}
\item[$(\alpha)$] $w^p\subseteq w^q$,\quad $u^p\subseteq u^q$, \qquad and
\item[$(\beta)$]  $(\forall (i,\xi)\in u^p)(f^p_{i,\xi}\subseteq
f^q_{i,\xi})$,\qquad and 
\item[$(\gamma)$] for each $(i,\xi)\in u^q$ one of the following occurs:

{\bf either}\ $f^q_{i,\xi}\rest u^p={\bf 0}_{u^p}$,

{\bf or}\ $i\in w^p$ and for some $\zeta,\vare<\chi_i$ we have $(i,\zeta)\in 
u^p$ and $f^q_{i,\xi}\rest u^p=(f^p_{i,\zeta})^\vare$, where $(f^p_{i,\zeta}
)^\vare:u^p\longrightarrow 2$ is defined by  
\[(f^p_{i,\zeta})^\vare(j,\gamma)=\left\{\begin{array}{ll}
0                       &\mbox{if } j=i,\ \vare\leq\gamma<\chi_i,\\
f^p_{i,\zeta}(j,\gamma) &\mbox{otherwise,}
\end{array}\right.\]

{\bf or}\ $i\notin w^p$ and {\em either} $f^q_{i,\xi}\rest
u^p=(f^p_{j,\zeta})^\vare$ (defined above) for some $(j,\zeta)\in u^p$,
$\vare<\chi_j$ {\em or} $f^q_{i,\xi}\rest u^p= (f^p_{j,\zeta})_{j'}$ for
some $(j,\zeta)\in u^p$ and $j'\leq j$, where $(f^p_{j,\zeta})_{j'}:u^p
\longrightarrow 2$ is defined by
\[(f^p_{j,\zeta})_{j'}(j^*,\gamma^*)=\left\{\begin{array}{ll}
0                           &\mbox{if } j'\leq j^*,\\
f^p_{j,\zeta}(j^*,\gamma^*) &\mbox{otherwise.}
\end{array}\right.\]
\end{enumerate}
\item Conditions $p,q\in\Ps$ are said to be {\em isomorphic} if
the well orderings 
\[(u^p,\prec_S\rest u^p)\quad\mbox{ and }\quad(u^q,\prec_S\rest u^q)\]
are isomorphic, and if $H:u^p\longrightarrow u^q$ is the
$\prec_S$--isomorphism then: 
\begin{enumerate}
\item[$(\alpha)$] $H(i,\xi)=(j,0)$ if and only if $\xi=0$,
\item[$(\beta)$]  $f^p_{i,\xi}=f^q_{H(i,\xi)}\comp H$ (for $(i,\xi)\in
u^p$). 
\end{enumerate}
\end{enumerate}
\end{definition}

\begin{proposition}
\label{2.3}
Let $S=(\mu,\lambda,\bar{\chi})$ be a good parameter. Then $\Ps$ is a
$(<\mu)$--complete $\mu^+$--cc forcing notion.
\end{proposition}

\begin{proof}
Plainly $\Ps$ is a $(<\mu)$--complete forcing notion (compare the proof of
\ref{1.6}). To verify the chain condition suppose that $\cA\subseteq\Ps$,
$|\cA|=\mu^+$. Apply the $\Delta$--lemma and ``standard cleaning'' to choose 
isomorphic conditions $p,q\in\cA$ such that if $H:u^p\longrightarrow u^q$ is 
the isomorphism from $p$ to $q$ then $H\rest (u^p\cap u^q)$ is the identity
on $u^p\cap u^q$. Put $w^r=w^p\cup w^q$, $u^r=u^p\cup u^q$ and for
$(i,\xi)\in u^r$ define a function $f^r_{i,\xi}: u^r\longrightarrow 2$ as
follows.
\begin{itemize}
\item If $(i,\xi)\in u^p$, $i\in w^p\cap w^q$ then $f^r_{i,\xi}=
f^p_{i,\xi}\cup (f^q_{H(i,\xi)})^{\xi+1}$,
\item if $(i,\xi)\in u^q$, $i\in w^p\cap w^q$ then $f^r_{i,\xi}=
(f^p_{H^{-1}(i,\xi)})^{\xi+1}\cup f^q_{i,\xi}$,
\item if $(i,\xi)\in u^p$, $i\in w^p\setminus w^q$ then $f^r_{i,\xi}=
f^p_{i,\xi}\cup (f^q_{H(i,\xi)})_i$, 
\item if $(i,\xi)\in u^q$, $i\in w^q\setminus w^p$ then $f^r_{i,\xi}=
(f^p_{H^{-1}(i,\xi)})_i\cup f^q_{i,\xi}$.
\end{itemize}
It is a routine to check that the functions $f^r_{i,\xi}$ are well defined 
and that they satisfy the demand \ref{2.2}(1c). Hence $r=\langle w^r,u^r,
\langle f^r_{i,\xi}: (i,\xi)\in u^r\rangle\rangle\in\Ps$ and one 
easily checks that it is an upper bound to both $p$ and $q$.  
\end{proof}

For a condition $p\in\Ps$ let 
\[F^p=\{(f^p_{i,\xi})^\vare,(f^p_{i,\xi})_j:(i,\xi)\in u^p,\ \vare<\chi_i,\
j\leq i\},\] 
where $(f^p_{i,\xi})^\vare,(f^p_{i,\xi})_j:u^p\longrightarrow 2$ are defined 
like in \ref{2.2}(1$\gamma$): 
\[(f^p_{i,\xi})^\vare(i',\zeta')=\left\{\begin{array}{ll}
0                       &\mbox{if } i=i',\ \vare\leq\zeta',\\
f^p_{i,\xi}(i',\zeta')&\mbox{otherwise,}
\end{array}\right.\]
\[(f^p_{i,\xi})_j(i',\zeta')=\left\{\begin{array}{ll}
0                           &\mbox{if } j\leq i',\\
f^p_{i,\xi}(i',\zeta') &\mbox{otherwise.}
\end{array}\right.\]
Like in the previous section, $\bB_p$ is the Boolean algebra $\bB_{(u^p,
F^p)}$ (see \ref{0.C}) (note that $p\leq q$ implies that $\bB_p$ is a
subalgebra of $\bB_q$). Let $\dot{\bB}_S^1$ be a $\Ps$--name such that
\[\forces_{\Ps}\mbox{`` }\dot{\bB}_S^1=\bigcup\{\bB_p: p\in\Gamma_{\Ps}\}
\mbox{ '',}\]
and for $s\in\xs$ let $\dot{f}_s$ be a $\Ps$--name such that 
\[\forces_{\Ps}\mbox{`` }\dot{f}_s=\bigcup\{f^p_s:s\in u^p,\ p\in
\Gamma_{\Ps}\}\mbox{ ''.}\] 

\begin{proposition}
\label{2.4}
Assume that $S=(\mu,\lambda,\bar{\chi})$ is a good parameter. Then in
$\V^{\Ps}$:  
\begin{enumerate}
\item For $s\in\xs$, $\dot{f}_s:\xs\longrightarrow 2$  is such that
$\dot{f}_s(s)=1$ and 
\[(\forall s'\in\xs)(s\prec_S s'\ \Rightarrow\ \dot{f}_s(s')=0).\] 
\item $\dot{\bB}_S^1$ is the Boolean algebra $\bB_{(\xs,\dot{F})}$ (see
\ref{0.C}), where   
\[\dot{F}=\{(\dot{f}_{i,\xi})^\vare, (\dot{f}_{i,\xi})_j:(i,\xi)\in\xs,\
\vare<\chi_i,\ j\leq i\},\] 
and $(\dot{f}_{i,\xi})^\vare, (\dot{f}_{i,\xi})_j:\xs\longrightarrow 2$ are 
such that 
\[(\dot{f}_{i,\xi})^\vare(i',\zeta')=\left\{\begin{array}{ll}
0                       &\mbox{if } i=i',\ \vare\leq\zeta',\\
\dot{f}_{i,\xi}(i',\zeta')&\mbox{otherwise,}
\end{array}\right.\]
\[(\dot{f}_{i,\xi})_j(i',\zeta')=\left\{\begin{array}{ll}
0                           &\mbox{if } j\leq i',\\
\dot{f}_{i,\xi}(i',\zeta') &\mbox{otherwise.}
\end{array}\right.\]
\item The sequence $\langle x_s: s\in\xs\rangle$ is left--separated in
$\dot{\bB}_S^1$ (when we consider $\xs$ with the well ordering $\prec_S$).  
\end{enumerate}
\end{proposition}

\begin{theorem}
\label{2.5}
Assume $S=(\mu,\lambda,\bar{\chi})$ is a good parameter. Then 
\[\forces_{\Ps}\mbox{`` there is no ideal $I\subseteq\dot{\bB}^1_S$ such
that $\pi(\dot{\bB}^1_S/I)=\lambda$ ''.}\]
\end{theorem}

\begin{proof}
Not surprisingly, the proof is similar to the one of \ref{1.8}. Let 
$\dot{I}$ be a $\Ps$--name for an ideal in $\dot{\bB}^1_S$, $p\in \Ps$, and 
suppose that $p\forces_{\Ps}\pi(\dot{\bB}^1_S/\dot{I})=\lambda$.   

Fix $i<\cf(\lambda)$. Use \ref{2.1} to choose $p_i$, $\theta_i$, $n_i$,
$\dcD_i$, $\dot{e}_i$ and $\dot{t}_i$ such that 
\begin{enumerate}
\item[$(\alpha)$] $p_i\in\Ps$ is a condition stronger than $p$, $\theta_i$
is a regular cardinal, $\chi_i^+<\theta_i<\lambda$ and $n_i\in\omega$,
\item[$(\beta)$]  $\dcD_i$ is a $\Ps$--name for a $(<\theta_i)$--complete
filter on $\theta_i$ extending the filter of co-bounded subsets of
$\theta_i$, 
\item[$(\gamma)$] $\forces_{\Ps}$`` $\dot{e}_i:\theta_i\times n_i
\longrightarrow\xs$ and $\dot{t}_i:\theta_i\times n_i\longrightarrow 2$ ''; 

for $\alpha<\theta_i$ let $\dot{a}^i_\alpha$ be a $\Ps$--name for an element 
of $\dot{\bB}^1_S$ such that 
\[\forces_{\Ps}\mbox{`` }\dot{a}^i_\alpha=\bigwedge\limits_{\ell<n_i} x^{
\dot{t}_i(\alpha,\ell)}_{\dot{e}_i(\alpha,\ell)}\mbox{ ''},\]
\item[$(\delta)$] $p_i\forces_{\Ps}$`` $\dot{a}^i_\alpha\in\dot{\bB}^1_S
\setminus \dot{I}$ ''\quad
for each $\alpha<\theta_i$, 
\item[$(\vare)$]  $p_i\forces_{\Ps}$`` if $b\in\dot{\bB}^1_S\setminus
\dot{I}$ then $\{\alpha<\theta_i: b\leq\dot{a}^i_\alpha\mod\dot{I}\}=
\emptyset\mod\dcD_i$\quad and

\qquad $(\forall\alpha<\theta_i)(\forall\beta<\alpha)(\dot{a}^i_\beta\wedge 
(-\dot{a}^i_\alpha)\notin\dot{I})$ ''. 
\end{enumerate}
For each $\alpha<\theta_i$ choose a maximal above $p_i$ antichain $\{p^i_{
\alpha,\zeta}:\zeta<\mu\}$ such that each $p^i_{\alpha,\zeta}\geq p_i$
decides the values of $\dot{e}_i(\alpha,\cdot)$, $\dot{t}_i(\alpha,\cdot)$. 
Let  
\[p^i_{\alpha,\zeta}\forces_{\Ps}\mbox{`` }\dot{e}_i(\alpha,\ell)=e^\zeta_i( 
\alpha,\ell)\quad\&\quad\dot{t}_i(\alpha,\ell)=t^\zeta_i(\alpha,\ell)\mbox{
''\qquad(for $\ell<n_i$),}\]
and we may assume that $(i,0), e^\zeta_i(\alpha,\ell)\in u^{p^i_{\alpha,
\zeta}}$ for $\alpha<\theta_i$, $\ell<n_i$ and $\zeta<\mu$. Take a generic
filter $G\subseteq\Ps$ such that $p_i\in G$ and work in $\V[G]$. Choose
ordinals $\dot{\gamma}^G_i<\theta_i$ and $\dot{\zeta}^G_i<\mu$ such that the 
set   
\[\begin{array}{ll}
\dot{X}^G_i\stackrel{\rm def}{=}\Big\{\beta<\theta_i:&\dot{\gamma}^G_i\leq
\beta\ \mbox{ and }\ p^i_{\dot{\gamma}^G_i,\dot{\zeta}^G_i}, p^i_{\beta,
\dot{\zeta}^G_i}\in G\ \mbox{ and }w^{p^i_{\dot{\gamma}^G_i,\dot{\zeta}^G_i
}}=w^{p^i_{\beta,\dot{\zeta}^G_i}},\\
&\mbox{the conditions }p^i_{\dot{\gamma}^G_i,\dot{\zeta}^G_i}, p^i_{\beta,
\dot{\zeta}^G_i}\mbox{ are isomorphic,\quad and}\\
&\mbox{if }H:u^{p^i_{\dot{\gamma}^G_i,\dot{\zeta}^G_i}}\longrightarrow u^{
p^i_{\beta,\dot{\zeta}^G_i}}\mbox{ is the isomorphism then}\\
&(\forall\ell<n_i)(H(e^{\dot{\zeta}^G_i}_i(\dot{\gamma}^G_i,\ell))=e^{\dot{
\zeta}^G_i}_i(\beta,\ell)\ \&\ t^{\dot{\zeta}^G_i}_i(\dot{\gamma}^G_i,\ell)= 
t^{\dot{\zeta}^G_i}_i(\beta,\ell))\\
&\mbox{and if }j\leq i,\ (j,\xi)\in\xs\ \mbox{ then}\\
&(j,\xi)\in u^{p^i_{\dot{\gamma}^G_i,\dot{\zeta}^G_i}}\ \Leftrightarrow\
(j,\xi)\in u^{p^i_{\beta,\dot{\zeta}^G_i}}\Big\} 
  \end{array}\]
is not $\emptyset$ modulo $\dcD^G_i$. Let $\dot{\delta}^G_i=\otp(u^{p^i_{
\dot{\gamma}^G_i,\dot{\zeta}^G_i}},\prec_S)$ and for $\alpha\in\dot{X}^G_i$
let $\langle s^{\alpha,i}_\vare:\vare<\dot{\delta}^G_i\rangle$ be the
$\prec_S$--increasing enumeration of $u^{p^i_{\alpha,\dot{\zeta}^G_i}}$. 
Apply Lemma \ref{0.A} to find a sequence $\langle s^{*,i}_\vare:\vare<
\dot{\delta}^G_i\rangle\subseteq\xs$ and a set $\dot{v}^G_i\subseteq
\dot{\delta}^G_i$ such that 
\begin{enumerate}
\item[(i)] $(\forall\vare\in\dot{\delta}^G_i\setminus\dot{v}^G_i)(\chi_i^+
\leq\cf(\{s\in\xs:s\prec_S s^{*,i}_\vare\},\prec_S)\leq\theta_i)$,
\item[(ii)] the set 
\[\begin{array}{ll}
\dot{B}^G_i\stackrel{\rm def}{=}\Big\{\beta\in\dot{X}^G_i:&\mbox{if }\vare
\in\dot{v}^G_i\mbox{ then } s^{\beta,i}_\vare=s^{*,i}_\vare,\ \ \mbox{
and}\\ 
&\mbox{if }\vare\in\dot{\delta}^G_i\setminus \dot{v}^G_i\mbox{ then}\\
&\sup_{\prec_S}\{s^{*,i}_\zeta:\zeta<\dot{\delta}^G_i,\ s^{*,i}_\zeta\prec_S 
s^{*,i}_\vare\}\prec_S s^{\beta,i}_\vare\prec_S s^{*,i}_\vare\Big\}
  \end{array}\]
is not $\emptyset$ modulo the filter $\dcD^G_i$,
\item[(iii)]\  if $s'_\vare\prec_S s^{*,i}_\vare$ for $\vare\in
\dot{\delta}^G_i\setminus \dot{v}^G_i$ then 
\[\{\beta\in\dot{B}^G_i:(\forall\vare\in\dot{\delta}^G_i\setminus\dot{v}^G_i
)(s'_\vare\prec_S s^{\beta,i}_\vare)\}\neq\emptyset\mod\dcD^G_i.\]  
\end{enumerate}
We may assume that $\dot{\gamma}^G_i\in\dot{B}^G_i$.

Now, in $\V$, we choose a condition $q_i\in\Ps$, ordinals $\gamma_i,\zeta_i,
\delta_i$, a set $v_i$ and a sequence $\langle s^{*,i}_\vare:\vare<\delta_i
\rangle\subseteq\xs$ such that $q_i\geq p^i_{\gamma_i,\zeta_i}$, and $q_i$
forces that these objects are as described above. If some condition stronger
than $q_i$ forces that $\alpha\in\dot{B}_i$, then we will use $\langle
s^{\alpha,i}_\vare:\vare<\delta_i\rangle$ to denote the
$\prec_S$--increasing enumeration of $u^{p^i_{\alpha,\zeta_i}}$.

Next, like in the proof of \ref{1.8}, we pick an unbounded set $Y\subseteq
\cf(\lambda)$ and $n<\omega$, $\delta<\mu$, $v\subseteq\delta$ such that for
$i,j\in Y$: 
\begin{itemize}
\item $n_i=n$, $\delta_i=\delta$, $v_i=v$, and
\item the conditions $p^i_{\gamma_i,\zeta_i},p^j_{\gamma_j,\zeta_j}$ are
isomorphic, and the isomorphism maps $e^{\zeta_i}_i(\gamma_i,\cdot)$ and
$t^{\zeta_i}_i(\gamma_i,\cdot)$ onto $e^{\zeta_j}_j(\gamma_j,\cdot),
t^{\zeta_j}_j(\gamma_j,\cdot)$, respectively.
\end{itemize}
Now use Lemma \ref{0.A} to find a sequence $\langle s_{*,\vare}:\vare<\delta
\rangle\subseteq\xs\cup\{(\cf(\lambda),0)\}$ and a set $v^*\subseteq\delta$
such that
\begin{enumerate}
\item[(a)] $(\forall\vare\in\delta\setminus v^*)(\cf(\{s\in\xs:s\prec_S
s_{*,\vare}\},\prec_S)=\cf(\lambda))$,
\item[(b)] the set 
\[\begin{array}{ll}
C\stackrel{\rm def}{=}\Big\{i\in Y:&\mbox{if }\vare\in v^*\mbox{ then }
s^{*,i}_\vare=s_{*,\vare},\ \ \mbox{ and}\\ 
&\mbox{if }\vare\in\delta\setminus v^*\mbox{ then}\\
&\sup_{\prec_S}\{s_{*,\zeta}:\zeta<\delta,\ s_{*,\zeta}\prec_S s_{*,\vare} 
\}\prec_S s^{*,i}_\vare\prec_S s_{*,\vare}\Big\}
  \end{array}\]
is unbounded in $\cf(\lambda)$,
\item[(c)] if $s'_\vare\prec_S s_{*,\vare}$ for $\vare\in\delta\setminus
v^*$, then the set  
\[\{i\in C:(\forall\vare\in\delta\setminus v^*)(s'_\vare\prec_S
s^{*,i}_\vare)\}\] 
is unbounded in $\cf(\lambda)$. 
\end{enumerate}
Next choose a set $C^+\in [C]^{\textstyle \cf(\lambda)}$ and ordinals
$\alpha_i<\beta_i<\theta_i$ (for $i\in C^+$) such that for every $i\in C^+$:
\begin{enumerate}
\item[(d)] if $\vare\in\delta\setminus v^*$ then for all $j\in C^+\cap i$
and $\zeta<\delta$ we have
\[\begin{array}{rcl}
s^{*,j}_\zeta\prec_S s_{*,\vare}&\Rightarrow& s^{*,j}_\zeta\prec_S
s^{*,i}_\vare,\quad\mbox{ and}\\
s^{\alpha_j,j}_\zeta\prec_S s_{*,\vare}&\Rightarrow& s^{\alpha_j,j}_\zeta
\prec_S s^{*,i}_\vare,\quad\mbox{ and}\\
s^{\beta_j,j}_\zeta\prec_S s_{*,\vare}&\Rightarrow& s^{\beta_j,j}_\zeta
\prec_S s^{*,i}_\vare,
  \end{array}\]
\item[(e)] some condition stronger than $q_i$ forces that $\alpha_i,\beta_i
\in\dot{B}_i$,  
\item[(f)] if $\vare\in\delta\setminus v$ and $x\in\{\alpha_i,\beta_i\}$,
then for all $j\in C^+\cap i$ and $\zeta<\delta$ we have 
\[\begin{array}{lcl}
s^{*,j}_\zeta\prec_S s^{*,i}_\vare\ \Rightarrow\ s^{*,j}_\zeta\prec_S
s^{x,i}_\vare&\mbox{ and }& s^{\alpha_j,j}_\zeta\prec_S s^{*,i}_\vare\
\Rightarrow\ s^{\alpha_j,j}_\zeta\prec_S s^{x,i}_\vare,\quad\mbox{and}\\
s^{\beta_j,j}_\zeta\prec_S s^{*,i}_\vare\ \Rightarrow\ s^{\beta_j,j}_\zeta
\prec_S s^{x,i}_\vare&\mbox{ and }& s^{\alpha_i,i}_\zeta\prec_S
s^{*,i}_\vare\ \Rightarrow\ s^{\alpha_i,i}_\zeta\prec_S s^{\beta_i,i}_\vare, 
  \end{array}\]
\item[(g)] if $\vare\in v^*$, $s_{*,\vare}=(j,\zeta)$ then $j<\min(C^+)$.
\end{enumerate}
Then $\langle\langle s^{\alpha_i,i}_\vare,s^{\beta_i,i}_\vare:\vare<\delta
\rangle:i\in C^+\rangle$ forms a $\Delta$--system of sequences with heart
$\langle s_{*,\vare}:\vare\in v\cap v^*\rangle$; but note that $s^{\alpha_i,
i}_\vare=s^{\beta_i,i}_\vare$ for $\vare\in v$. Let $u^*=\{s_{*,\vare}:\vare
\in v\cap v^*\}$ and $w^*=\{j<\cf(\lambda): (j,0)\in u^*\}$.

\begin{claim}
\label{2.5.1}
For each $i_0\in C^+$:
\[q_{i_0}\forces_{\Ps}\mbox{`` }(\forall\alpha\in\dot{B}_{i_0})(\exists i^*
\in C^+)(\dot{a}^{i^*}_{\alpha_{i^*}}\wedge(-\dot{a}^{i^*}_{\beta_{i^*}})
\leq\dot{a}^{i_0}_\alpha\ \&\ p_{i^*}\in\Gamma_{\Ps})\mbox{ ''}\] 
(where $\dot{B}_{i_0}$ was defined in {\rm (ii)}).  
\end{claim}

\begin{proof}[Proof of the claim]
Let $i_0\in C^+$. We will show that for every condition $q\geq q_{i_0}$ and
an ordinal $\alpha<\theta_{i_0}$ such that $q\forces\alpha\in\dot{B}_{i_0}$,
there are $i^*\in C^+$ and a condition $r$ stronger than both $q$ and
$p_{i^*}$, and such that $r\forces$`` $\dot{a}^{i^*}_{\alpha_{i^*}}\wedge
(-\dot{a}^{i^*}_{\beta_{i^*}})\leq\dot{a}^{i_0}_\alpha$ ''. 

So suppose $q\geq q_i$, $q\forces\alpha\in\dot{B}_{i_0}$. We may assume that 
$p^{i_0}_{\alpha,\zeta_{i_0}}\leq q$. Choose $i^*\in C^+\setminus (i_0+1)$
such that   
\[u^q\cap u^{p^{i^*}_{\alpha_{i^*},\zeta_{i^*}}}=u^q\cap
u^{p^{i^*}_{\beta_{i^*},\zeta_{i^*}}}=u^*\qquad\mbox{and}\qquad w^q\subseteq
i^*,\] 
Let $p^{i_0}_{\alpha,\zeta_{i_0}}=q^0$, $p^{i^*}_{\alpha_{i^*},\zeta_{i^*}}=
q^1$, $p^{i^*}_{\beta_{i^*},\zeta_{i^*}}=q^2$, and  
\[\tau_0=\bigwedge_{\ell<n}x^{t^{\zeta_{i_0}}_{i_0}(\alpha,\ell)}_{e^{
\zeta_{i_0}}_{i_0}(\alpha,\ell)},\qquad\tau_1=\bigwedge_{\ell<n}x^{t^{
\zeta_{i^*}}_{i^*}(\alpha_{i^*},\ell)}_{e^{\zeta_{i^*}}_{i^*}(\alpha_{
i^*},\ell)}\quad\mbox{ and }\quad\tau_2=\bigwedge_{\ell<n}x^{t^{\zeta_{
i^*}}_{i^*}(\beta_{i^*},\ell)}_{e^{\zeta_{i^*}}_{i^*}(\beta_{i^*},\ell)}\] 
(so $q^0\leq q$ and $\tau_0\in \bB_{q^0}\subseteq\bB_q$, $\tau_1\in
\bB_{q^1}$, $\tau_2\in\bB_{q^2}$). Note that the conditions $q^0,q^1,q^2$
are pairwise isomorphic and the isomorphism $H^{k,k'}$ from $q^k$ to
$q^{k'}$ carries $\tau_k$ to $\tau_{k'}$. Moreover, $H^{k,k'}$ is the
identity on $u^{q^k}\cap u^{q^{k'}}$. Also note that $w^{q^1}=w^{p^{i^*}_{
\gamma_{i^*},\zeta_{i^*}}}=w^{q^2}$ and, as $w^q\subseteq i^*$, our choices
imply $H^{k,0}(i,\xi)\preceq_S(i,\xi)$ for $k=1,2$, $(i,\xi)\in u^{q^k}$. 

Now we define a condition $r$ stronger than $q,q^1,q^2$. We put $w^r=w^q\cup 
w^{q^1}$, $u^r=u^q\cup u^{q^1}\cup u^{q^2}$ and we define functions $f^r_{i,
\xi}:u^r\longrightarrow 2$ as follows.

\begin{enumerate}
\item If $(i,\xi)\in u^{q^1}\cap u^{q^2}$, $i\in w^q$ then we let $f^r_{i,
\xi}=f^q_{H^{1,0}(i,\xi)}\cup f^{q^1}_{i,\xi}\cup f^{q^2}_{i,\xi}$. 

\noindent [Note that by (d)+(ii) we have $(i,0)\preceq_S H^{1,0}(i,\xi)
\preceq_S (i,\xi)$.] 
\item If $(i,\xi)\in u^{q^1}\cap u^{q^2}$, $i\notin w^q$ then we first
choose $\vare^*$ such that, if possible, $(f^{q^0}_{H^{1,0}(i,\xi)})^{
\vare^*}(\tau_0)=1$, and then we let $f^r_{i,\xi}=(f^q_{H^{1,0}(i,\xi)})^{
\vare^*}\cup f^{q^1}_{i,\xi}\cup f^{q^2}_{i,\xi}$. 

\noindent [Note that $H^{1,0}(i,\xi)\prec_S (i,\xi)$, and thus if $H^{1,0}
(i,\xi)=(j,\zeta)$ then $j<i$, $j\notin w^{q^1}$.] 
\item If $(i,\xi)\in u^{q^2}\setminus u^{q^1}$ (so $i>i^*\geq\sup(w^q)$)
then we first choose $\vare^*$ such that, if possible, $(f^{q^0}_{H^{2,0}
(i,\xi)})^{\vare^*}(\tau_0)=1$, and then we let $f^r_{i,\xi}=(f^q_{H^{2,0}
(i,\xi)})^{\vare^*}\cup f^{q^1}_{H^{2,1}(i,\xi)}\cup f^{q^2}_{i,\xi}$.

\noindent [Note that $H^{2,0}(i,\xi)\prec_S (i,0)\prec_S H^{2,1}(i,\xi)
\prec_S (i,\xi)$; remember $w^{q^1}=w^{q^2}$. Also, if $H^{2,0}(i,\xi)=(j,
\zeta)$, then $j\notin w^{q^1}$.] 
\item If $(i,\xi)\in u^{q^1}\setminus u^{q^2}$ then, like above, we choose
$\vare^*$ such that if possible then $(f^{q^0}_{H^{1,0}(i,\xi)})^{\vare^*}
(\tau_0)=1$, and next we put $f^r_{i,\xi}=(f^q_{H^{1,0}(i,\xi)})^{\vare^*}
\cup f^{q^1}_{i,\xi}\cup (f^{q^2}_{H^{1,2}(i,\xi)})^{\xi+1}$. 
\item If $(i,\xi)\in u^q\setminus u^{q^1}$ then we look at $f^q_{i,\xi}\rest 
u^{q^0}$. If it is ${\bf 0}_{u^{q^0}}$ then we let $f^r_{i,\xi}=f^q_{i,\xi}
\cup {\bf 0}_{u^{q^1}}\cup {\bf 0}_{u^{q^2}}$. Otherwise, we consider the
following three cases.
\begin{enumerate}
\item[$(\alpha)$] Suppose $i\in w^{q^0}$. Then for some $\vare\leq\zeta<
\chi_i$, $\vare\leq\xi+1$ we have $f^q_{i,\xi}\rest u^{q^0}=(f^{q^0}_{i,
\zeta})^\vare$ and we let:

-- if $i\in w^{q^1}$ then $f^r_{i,\xi}=f^q_{i,\xi}\cup (f^{q^1}_{H^{0,1}
(i,\zeta)})^\vare\cup (f^{q^2}_{H^{0,2}(i,\zeta)})^{\vare}$,

-- if $i\notin w^{q^1}$ then $f^r_{i,\xi}=f^q_{i,\xi}\cup (f^{q^1}_{H^{0,1}
(i,\zeta)})_i \cup (f^{q^2}_{H^{0,2}(i,\zeta)})_i$. 
\end{enumerate}
[Note that if $i\in w^{q^1}$ then $(i,\zeta)\preceq_S H^{0,1}(i,\zeta)=
H^{0,2}(i,\zeta)\prec_S (i+1,0)$, and if $i\notin w^{q^1}$ then $(j,0)
\preceq_S H^{0,1}(i,\zeta)\preceq_S H^{0,2}(i,\zeta)\prec_S(j+1,0)$ for some 
$j>i$.] 

\begin{enumerate}
\item[$(\beta)$] Suppose $i\notin w^{q^0}$ (so $i\notin w^{q^1}$) and
$f^q_{i,\xi}\rest u^{q^0}=(f^{q^0}_{i',\zeta'})^{\vare'}$, $(i',\zeta')\in
u^{q^0}$, $\vare'\leq\zeta'<\chi_{i'}$.

-- If $i'\in w^{q^1}$ and $i'<i$, then put $f^r_{i,\xi}=f^q_{i,\xi}\cup
(f^{q^1}_{H^{0,1}(i',\zeta')})^{\vare'}\cup (f^{q^2}_{H^{0,2}(i',
\zeta')})^{\vare'}$. 

-- If $i'\in w^{q^1}$ and $i<i'$, then we put $f^r_{i,\xi}=f^q_{i,\xi}\cup
(f^{q^1}_{H^{0,1}(i',\zeta')})_i\cup (f^{q^2}_{H^{0,2}(i',\zeta')})_i$.  

-- If $i'\notin w^{q^1}$, then let $f^r_{i,\xi}=f^q_{i,\xi}\cup
(f^{q^1}_{H^{0,1}(i',\zeta')})_i \cup (f^{q^2}_{H^{0,2}(i',\zeta')})_i$.
\item[$(\gamma)$] Suppose $i\notin w^{q^0}$ and $f^q_{i,\xi}\rest u^{q^0}= 
(f^{q^0}_{i',\zeta'})_{j'}$, $j'\leq \min\{i,i'\}$, $(i',\zeta')\in
u^{q^0}$. Let $f^r_{i,\xi}=f^q_{i,\xi}\cup (f^{q^1}_{H^{0,1}(i',
\zeta')})_{j'}\cup (f^{q^2}_{H^{0,2}(i',\zeta')})_{j'}$.
\end{enumerate}
\end{enumerate}
Verifying that the functions $f^r_{i,\xi}$ are well defined and that $r= 
\langle w^r,u^r,\langle f^r_{i,\xi}: (i,\xi)\in u^r\rangle\rangle\in\Ps$ is
a condition stronger than $q,q^1,q^2$ is left to the reader. Let us argue
that $\bB_r\models \tau_1\wedge(-\tau_2)\leq \tau_0$. If not then we have a 
function $f\in F^r$ such that $f(\tau_0)=f(\tau_2)=0$ and $f(\tau_1)=1$. 
Clearly $f$ cannot be ${\bf 0}_{u^r}$, so it is either $(f^r_{i,\xi})^\vare$
or $(f^r_{i,\xi})_j$. Let us look at the definition of the functions
$f^r_{i,\xi}$ and consider each case there separately.  
\medskip

\noindent{\sc Cases 1, 5$\alpha,\beta,\gamma$:}\quad Plainly
$f^r_{i,\xi}(\tau_1)=f^r_{i,\xi}(\tau_2)$ and also $(f^r_{i,\xi})_j(\tau_1)
=(f^r_{i,\xi})_j(\tau_2)$ (remember $w^{q^1}=w^{q^2}$). As far as the
operation $(\cdot)^\vare$ is concerned, note that $(\{i\}\times\chi_i)\cap
u^{q^1}=(\{i\}\times\chi_i)\cap u^{q^2}$, so (in these cases) we easily get
$(f^r_{i,\xi})^\vare (\tau_1)=(f^r_{i,\xi})^\vare(\tau_2)$, a
contradiction. 
\medskip

\noindent{\sc Case 2:}\quad Again, $f^r_{i,\xi}(\tau_1)=f^r_{i,\xi}(\tau_2)$ 
and $(f^r_{i,\xi})_j(\tau_1)=(f^r_{i,\xi})_j(\tau_2)$ (for each $j$). So
suppose that $f=(f^r_{i,\xi})^\vare$ for some $\vare$, and look at the
choice of $\vare^*$ in the current case. Since $1=(f^r_{i,\xi})^\vare(
\tau_1)=(f^{q^1}_{i,\xi})^\vare(\tau_1)$, we conclude that $1=(f^{q^0}_{
H^{1,0}(i,\xi)})^{\vare^*}(\tau_0)=f^r_{i,\xi}(\tau_0)=(f^r_{i,\xi})^\vare
(\tau_0)$, a contradiction. 
\medskip

\noindent{\sc Case 3:}\quad Note that $f^r_{i,\xi}(\tau_1)=f^r_{i,\xi}(
\tau_2)$ (and also $(f^r_{i,\xi})_j(\tau_1)=(f^r_{i,\xi})_j(\tau_2)$). Now, 
if for some $\vare$ we have $(f^r_{i,\xi})^\vare(\tau_1)=1$, then look at
the choice of $\vare^*$ -- necessarily $(f^r_{i,\xi})^\vare(\tau_0)=
f^r_{i,\xi}(\tau_0)=1$ (remember $(i,0)\prec_S H^{2,1}(i,\xi)\prec_S
(i,\xi)$). 
\medskip

\noindent{\sc Case 4:}\quad Like above: if for some $\vare$ we have
$(f^r_{i,\xi})^\vare(\tau_1)=1$, then necessarily $f^r_{i,\xi}(\tau_0)=
(f^r_{i,\xi})^\vare(\tau_0)=1$. Moreover, $(f^r_{i,\xi})_j(\tau_1)=
(f^r_{i,\xi})_j(\tau_2)$ for all $j\leq i$. 
\medskip

In all cases we get a contradiction showing that $\bB_r\models\tau_1\wedge
(-\tau_2)\leq \tau_0$, and hence $r\forces$`` $\dot{a}^{i^*}_{\alpha_{i^*}}
\wedge(-\dot{a}^{i^*}_{\beta_{i^*}})\leq\dot{a}^{i_0}_\alpha$ '', finishing
the proof of the claim.  
\end{proof}

Finally we note that \ref{2.5.1} and clauses $(\beta)$, $(\vare)$ give an
immediate contradiction, showing the theorem. 
\end{proof}

\begin{conclusion}
\label{2.6}
It is consistent that there is a Boolean algebra $\bB$ of size $\lambda$
such that there is a left--separated sequence of length $\lambda$ in $\bB$
(and thus $\hd^+_{(5)}(\bB)=\lambda^+$), but there is no ideal $I\subseteq
\bB$ with $\pi(\bB/I)=\lambda$ (so $\hd^+_{(7)}(\bB)=\hd_{(7)}(\bB)=
\lambda$). 
\end{conclusion}

\begin{problem}
Can one construct a Boolean algebra $\bB$ as in \ref{2.6} for $\lambda$ from
any cardinal arithmetic assumptions?
\end{problem}

\section{More on the attainment problem}

In this section we will assume the following:
\begin{hypothesis}
\label{Hyp3.1}
$S=(\mu,\lambda,\bar{\chi})$ is such that $\mu,\lambda$ are cardinals
satisfying 
\[\mu=\mu^{<\mu}<\cf(\lambda)<\lambda\leq 2^\mu,\]
and $\bar{\chi}=\langle\chi_i:i<\cf(\lambda)\rangle$ is a strictly
increasing continuous sequence of cardinals such that 
\[\chi_0=0,\quad \cf(\lambda)<\chi_1,\quad \cf(\chi_{i+1})=\chi_{i+1},\quad
\mbox{ and }\quad \sup_{i<\cf(\lambda)}\chi_i=\lambda.\]
For $\alpha<\lambda$ let $j(\alpha)<\cf(\lambda)$ be such that
$\chi_{j(\alpha)}\leq \alpha<\chi_{j(\alpha)+1}$.
\end{hypothesis}

\begin{definition}
\label{Def3.2}
\begin{enumerate}
\item A pair $(\bar{\eta},A)$ is {\em a base for $S=(\mu,\lambda,
\bar{\chi})$} if 
\begin{enumerate}
\item[(a)] $A\subseteq\mu^{\textstyle{<}\mu}$, $\bar{\eta}=\langle
\eta_\alpha:\alpha<\lambda\rangle\subseteq\mu^{\textstyle \mu}$,
\item[(b)] if $\alpha<\beta<\lambda$, $j(\alpha)=j(\beta)$ then
$\eta_\alpha\cap\eta_\beta\notin A$, and
\item[(c)] if $Y\in [\lambda]^{\textstyle\lambda}$ then there are distinct
$\alpha,\beta\in Y$ such that $\eta_\alpha\cap\eta_\beta\in A$.
\end{enumerate}
\item $(\bar{\eta},A)$ is called {\em a base$^+$ for $S$} if it satisfies
demands (a), (b) (stated above) and
\begin{enumerate}
\item[(c$^+$)] if $Y\in [\lambda]^{\textstyle\lambda}$ and ${\bf t}\in
\{0,1\}$, then there are $\alpha,\beta\in Y$ such that 
\[\alpha<\beta,\qquad\eta_\alpha\cap\eta_\beta\in A,\qquad\mbox{ and }\qquad 
\eta_\alpha<_\lex\eta_\beta\ \mbox{ iff }\ {\bf t}=0.\]
\end{enumerate}
\end{enumerate}
\end{definition}

For a topological space $X$, a $(\kappa_0,\kappa_1)$--Lusin set in $X$ is a
set $L\subseteq X$ such that $|L|=\kappa_0$ and for every meager subset $Z$
of  $X$ the intersection $Z\cap L$ is of size less than $\kappa_1$. (See,
e.g., Cicho\'n \cite{Ci89} for a discussion of sets of this type.) Below,
the space $\mu^{\textstyle\mu}$ is equipped with the topology generated by
sets of the form 
\[[\rho]=\{\eta\in\mu^{\textstyle\mu}:\rho\vartriangleleft\eta\}\]
for $\rho\in \mu^{\textstyle{<}\mu}$.

\begin{proposition}
\label{Prop3.3}
Assume that for some $i^*<\cf(\lambda)$ there is a $(\lambda,
\chi_{i^*})$--Lusin set $L$ in $\mu^{\textstyle\mu}$. Then there is a
base$^+$ for $S$.
\end{proposition}

\begin{proof}
Choose sequences $\langle\nu_i:i<\cf(\lambda)\rangle\subseteq\mu^{\textstyle
\mu}$ and $\langle \rho_\alpha:\alpha<\lambda\rangle\subseteq L$, both with
no repetitions. For $\alpha<\lambda$ let $\eta_\alpha\in\mu^{\textstyle\mu}$ 
be defined by
\[\eta_\alpha(2\cdot\xi)=\nu_{j(\alpha)}(\xi)\quad\mbox{ and }\quad
\eta_\alpha (2\cdot\xi+1)=\rho_\alpha(\xi)\] 
(for $\xi<\mu$), and let $A=\bigcup\limits_{\xi<\mu}\mu^{\textstyle 2\cdot
\xi}$. We claim that $(\langle\eta_\alpha:\alpha<\lambda\rangle,A)$ is a
base$^+$ for $S$. The conditions \ref{Def3.2}(1)(a,b) should be clear. Let
us verify \ref{Def3.2}(2)(c$^+$). So suppose that $Y\in
[\lambda]^{\textstyle \lambda}$ and ${\bf t}\in\{0,1\}$. Choose sequences
$\langle Y_i:i<\cf(\lambda)\rangle$ and $\langle j_i:i<\cf(\lambda)\rangle$
such that 
\begin{itemize}
\item $Y_i\subseteq Y$, $(\forall\alpha\in Y_i)(j(\alpha)=j_i)$, and
$|Y_i|=\chi_{i^*}$ (so $\{\rho_\alpha:\alpha\in Y_i\}$ is not meager),
\item the sequence $\langle j_i:i<\cf(\lambda)\rangle$ is strictly
increasing. 
\end{itemize}
For each $i<\cf(\lambda)$ pick $\sigma_i\in\mu^{\textstyle{<}\mu}$ such that 
\[(\forall\sigma\in \mu^{\textstyle{<}\mu})(\sigma_i\vartriangleleft\sigma
\ \ \Rightarrow\ \ [\sigma]\cap\{\rho_\alpha:\alpha\in Y_i\}\neq
\emptyset).\]
We may pick $i_0<i_1<\cf(\lambda)$ such that 
\[\sigma_{i_0}=\sigma_{i_1}=\sigma^*\qquad\mbox{ and }\qquad\nu_{j_{i_0}}
<_\lex \nu_{j_{i_1}}\ \ \mbox{ iff }\ \ {\bf t}=0.\] 
(Remember that, under the assumptions of \ref{Hyp3.1}, $(\mu^{\textstyle
\mu},<_\lex)$ contains no monotonic sequences of length $\cf(\lambda)$.)
Let $\xi=\lh(\nu_{j_{i_0}}\cap\nu_{j_{i_1}})$ and take $\sigma'\in\mu^{
\textstyle{<}\mu}$ such that $\sigma^*\trianglelefteq\sigma'$ and $\xi<
\lh(\sigma')$. Now pick $\alpha_0\in Y_{i_0}$ and $\alpha_1\in Y_{i_1}$ such
that $\sigma'\vartriangleleft\rho_{\alpha_0}\cap\rho_{\alpha_1}$ (there are
such $\alpha_0,\alpha_1$ by the choice of $\sigma_{i_0}=\sigma_{i_1}=
\sigma^*$). Note that then necessarily $\alpha_0<\alpha_1$, $\lh(
\eta_{\alpha_0}\cap\eta_{\alpha_1})=2\cdot\xi$ (so $\eta_{\alpha_0}\cap
\eta_{\alpha_1}\in A$) and $\eta_{\alpha_0}<_\lex\eta_{\alpha_1}$ iff ${\bf
t}=0$.  
\end{proof}

\begin{proposition}
\label{Prop3.4}
Let $\bP=(2^{\textstyle{<}\mu},\vartriangleleft)$ be the $\mu$--Cohen
forcing notion. Then
\[\forces_{\bP}\mbox{`` there is a base$^+$ for $S$ (and $S$ is still as in
\ref{Hyp3.1}) ''.}\]
\end{proposition}

\begin{proof}
Pick sequences $\langle\nu_i:i<\cf(\lambda)\rangle$ and $\langle\rho_\alpha:
\alpha<\lambda\rangle$ of pairwise distinct elements of $\mu^{\textstyle
\mu}$. Let $\dot{A}^*$ be a $\bP$--name for the generic subset of $\mu$
(added by $\bP$) and let $\dot{A}$ be a $\bP$--name such that
\[\forces_\bP\mbox{`` }\dot{A}=\{\nu\in\mu^{\textstyle{<}\mu}:\lh(\nu)\in
\dot{A}^*\}\mbox{ ''.}\]
For $\alpha<\lambda$, let $\dot{\eta}_\alpha$ be a $\bP$--name for a
function in $\mu^{\textstyle\mu}$ such that 
\[\forces_\bP(\forall\xi\!\in\!\dot{A}^*)(\dot{\eta}_\alpha(\xi)
=\nu_{j(\alpha)}(\otp(\dot{A}^*\cap\xi))\ \&\ (\forall\xi\!\in\!\mu\!
\setminus\!\dot{A}^*)(\dot{\eta}_\alpha(\xi)=\rho_\alpha(\otp(\xi\setminus
\dot{A}^*))).\] 
We claim that
\[\forces_\bP\mbox{`` }(\langle\dot{\eta}_\alpha:\alpha<\lambda\rangle,A)
\mbox{ is a base$^+$ for $S$ ''.}\]
Clauses \ref{Def3.2}(1)(a,b) should be clear, so let us prove
\ref{Def3.2}(2)(c$^+$) only. Let $\langle\dot{\alpha}_\gamma:\gamma<\lambda
\rangle$ be a $\bP$--name for an increasing $\lambda$--sequence of elements
of $\lambda$, and let ${\bf t}\in\{0,1\}$, $p\in\bP$. For each $\gamma<
\lambda$ pick a condition $p_\gamma\geq p$ and an ordinal $\alpha_\gamma$
such that $p_\gamma\forces\dot{\alpha}_\gamma=\alpha_\gamma$. Necessarily,
there are $X\in [\lambda]^{\textstyle\lambda}$ and $p^*\in\bP$ such that
$p^*=p_\gamma$ for $\gamma\in X$. Then also $\alpha_{\gamma_0}<
\alpha_{\gamma_1}$ for $\gamma_0<\gamma_1$ from $X$. Shrinking $X$ a little
we may also demand that for some sequences $\sigma_j\in\mu^{\textstyle
\lh(p^*)+2}$ (for $j<\cf(\lambda)$) we have
\[\gamma\in X\ \&\ j(\alpha_\gamma)=j\quad\Rightarrow\quad\sigma_j
\vartriangleleft \rho_{\alpha_\gamma}.\]
Now pick $\gamma_0<\gamma_1$ from $X$ such that letting $j_0=j(
\alpha_{\gamma_0})$ and $j_1=j(\alpha_{\gamma_1})$ we have
\[j_0<j_1\quad \mbox{ and }\quad \sigma_{j_0}=\sigma_{j_1}\quad \mbox{
and }\quad \nu_{j_0}<_\lex\nu_{j_1}\ \mbox{ iff }\ {\bf t}=0.\]
Let a condition $q\geq p^*$ be such that $\lh(q)=\lh(p^*)+\lh(\nu_{j_0}\cap 
\nu_{j_1})+2$ and $q(\xi)=1$ for all $\xi\in\lh(q)\setminus\lh(p^*)$. It
should be clear that $\alpha_{\gamma_0}<\alpha_{\gamma_1}$ and 
\[q\forces\mbox{`` }\dot{\eta}_{\alpha_{\gamma_0}}\cap\dot{\eta}_{\alpha_{
\gamma_1}}\in\dot{A}\quad\mbox{ and }\quad\dot{\eta}_{\alpha_{\gamma_0}}
<_\lex\dot{\eta}_{\alpha_{\gamma_1}}\ \mbox{ iff }\ {\bf t}=0\mbox{ ''.}\] 
\end{proof}

\begin{definition}
\label{BAbase}
Let ${\bf b}=(\bar{\eta},A)$ be a base for $S$, $\bar{\eta}=\langle
\eta_\alpha:\alpha<\lambda\rangle$. We define the Boolean algebra $\bB^{\bf
b}$ determined by ${\bf b}$. First, functions $f^{\bf b}_\alpha:\lambda
\longrightarrow 2$ (for $\alpha<\lambda$) are such that
\[f^{\bf b}_\alpha(\beta)=\left\{\begin{array}{ll}
1&\mbox{ if }\alpha=\beta\ \mbox{ or }\ \alpha\neq\beta\ \&\ \eta_\alpha\cap
\eta_\beta\in A\ \&\ \eta_\alpha<_\lex\eta_\beta,\\
0&\mbox{ otherwise.}
				 \end{array}\right.\]
Next, we let $F^{\bf b}=\{f^{\bf b}_\alpha:\alpha<\lambda\}$ and $\bB^{\bf
b}=\bB_{(\lambda,F^{\bf b})}$ (see \ref{0.C}).
\end{definition}

\begin{theorem}
\label{Thm3.5}
If ${\bf b}$ is a base for $S=(\mu,\lambda,\bar{\chi})$, then
\[\hL(\bB^{\bf b})=\hd(\bB^{\bf b})=s^+(\bB^{\bf b})=\lambda.\]
If additionally ${\bf b}$ is a base$^+$ for $S$ then also
\[\hL_{(7)}^+(\bB^{\bf b})=\hd_{(5)}^+(\bB^{\bf b})=\lambda.\]
\end{theorem}

\begin{proof}
Let ${\bf b}=(\bar{\eta},A)$, $\bar{\eta}=\langle\eta_\alpha:\alpha<\lambda
\rangle$. Clearly $|\bB^{\bf b}|=\lambda$.

\begin{claim}
\label{clx1}
$\hL(\bB^{\bf b})=\hd(\bB^{\bf b})=s(\bB^{\bf b})=\lambda$.
\end{claim}

\begin{proof}[Proof of the claim]
By \ref{Def3.2}(1)(b), $f_\alpha^{\bf b}(\beta)=0$ whenever
$\alpha\neq\beta$ and $j(\alpha)=j(\beta)$. Therefore, by \ref{xxx}(1), the
sequence $\langle x_\alpha:\chi_i\leq\alpha<\chi_{i+1}\rangle$ is ideal
independent (for each $i<\cf(\lambda)$).
\end{proof}

The main part is to show that $s^+(\bB^{\bf b})=\lambda$ (and/or under the
additional assumption, that $\hL_{(7)}^+(\bB^{\bf b})=\hd_{(5)}^+(\bB^{\bf
b})=\lambda$), and for this we will need the following technical claim. 

\begin{claim}
\label{clx2}
Suppose that $k^*,\ell^*<\omega$, $\alpha_k,\alpha_{\ell,k}<\lambda$ (for
$k<k^*$, $\ell<\ell^*$) and $\sigma_0,\ldots,\sigma_{k^*-1}\in\mu^{
\textstyle {<}\mu}$ are such that
\begin{enumerate}
\item[$(\alpha)$] $\sigma_0,\ldots,\sigma_{k^*-1}$ are pairwise
incomparable, 
\item[$(\beta)$] $\sigma_k\vartriangleleft\eta_{\alpha_k}$,
$\sigma_k\vartriangleleft\eta_{\alpha_{\ell,k}}$ (for $\ell<\ell^*$,
$k<k^*$),
\item[$(\gamma)$] for each $k<k^*$ one of the following occurs:
\begin{enumerate}
\item[(i)]  $\alpha_k=\alpha_{\ell,k}$ for some $\ell<\ell^*$, or
\item[(ii)] there are $\ell_1,\ell_2,\ell_3<\ell^*$ such that 
\begin{itemize}
\item $\eta_{\alpha_k}\cap\eta_{\alpha_{\ell_1,k}}\vartriangleleft
\eta_{\alpha_k}\cap\eta_{\alpha_{\ell_2,k}}\vartriangleleft
\eta_{\alpha_k}\cap\eta_{\alpha_{\ell_3,k}}$, and 
\item $\eta_{\alpha_k}\cap\eta_{\alpha_{\ell_1,k}},\eta_{\alpha_k}\cap
\eta_{\alpha_{\ell_2,k}}\in A$, and 
\item $\eta_{\alpha_{\ell_1,k}}<_\lex\eta_{\alpha_k}<_\lex\eta_{\alpha_{
\ell_2,k}}$.
\end{itemize}
\end{enumerate}
\end{enumerate}
Let $t(k)\in\{0,1\}$ for $k<k^*$. Then
\[\bB^{\bf b}\models\bigwedge_{k<k^*} x^{t(k)}_{\alpha_k}\leq
\bigvee_{\ell<\ell^*}\bigwedge_{k<k^*} x^{t(k)}_{\alpha_{\ell,k}}.\]
\end{claim}

\begin{proof}[Proof of the claim]
We are going to show that, under our assumptions, for each $f\in F^{\bf b}$
there is $\ell<\ell^*$ such that $(\forall k<k^*)(f(\alpha_k)=f(
\alpha_{\ell,k}))$. So let us fix $\beta<\lambda$, and we consider $f^{\bf
b}_\beta$. First note that 
\begin{enumerate}
\item[$(\boxtimes_k)$] if $\sigma_k$ is not an initial segment of
$\eta_\beta$, then $f^{\bf b}_\beta(\alpha_k)=f^{\bf b}_\beta(\alpha_{\ell,
k})$ for all $\ell<\ell^*$.
\end{enumerate}
[Why? Suppose $\sigma_k\ntriangleleft\eta_\beta$. Then clearly $\alpha_k\neq
\beta\neq\alpha_{\ell,k}$ (for $\ell<\ell^*$) and 
\[\eta_{\alpha_k}\cap\eta_\beta=\eta_{\alpha_{\ell,k}}\cap\eta_\beta\quad
\mbox{ and }\quad\eta_{\alpha_k}<_\lex\eta_\beta\ \Leftrightarrow\
\eta_{\alpha_{\ell,k}}<_\lex\eta_\beta.\]
Now look at the definition of $f^{\bf b}_\beta$.]\\
If no $\sigma_k$ is an initial segment of $\eta_\beta$, then (by
$(\boxtimes_k)$) we conclude $f^{\bf b}_\beta(\alpha_k)=f^{\bf b}_\beta(
\alpha_{\ell,k})$ for all $\ell<\ell^*$, $k<k^*$. So suppose that $\sigma_m
\vartriangleleft\eta_\beta$, $m<k^*$. Then for all $k<k^*$, $k\neq m$, we
have $\sigma_k\ntriangleleft\eta_\beta$ and thus $f^{\bf b}_\beta(\alpha_k)
=f^{\bf b}_\beta(\alpha_{\ell,k})$ (for all $\ell<\ell^*$). Thus it is
enough to find $\ell<\ell^*$ such that $f^{\bf b}_\beta(\alpha_m)= f^{\bf
b}_\beta(\alpha_{\ell,m})$. If $\alpha_m=\alpha_{\ell,m}$ for some
$\ell<\ell^*$, then this $\ell$ works. So suppose $\alpha_m\neq\alpha_{\ell,
m}$ for all $\ell<\ell^*$. Then clause $(\gamma)$(ii) holds true for $m$,
and let $\ell_1,\ell_2,\ell_3$ be as there. If $\eta_{\alpha_m}\cap
\eta_\beta\vartriangleleft\eta_{\alpha_m}\cap\eta_{\alpha_{\ell_3,m}}$, then
easily $f^{\bf b}_\beta(\alpha_m)=f^{\bf b}_\beta(\alpha_{\ell_3,m})$. 
Otherwise $\eta_{\alpha_m}\cap\eta_{\alpha_{\ell_3,m}}\trianglelefteq
\eta_{\alpha_m}\cap\eta_\beta$, and $f^{\bf b}_\beta(\alpha_{\ell_1,m})\neq
f^{\bf b}_\beta(\alpha_{\ell_2,m})$, so either $\ell_1$ or $\ell_2$ works.
\end{proof}

\begin{claim}
\label{clx3}
$s^+(\bB^{\bf b})=\lambda$.
\end{claim}

\begin{proof}[Proof of the claim]
Suppose that $\langle a_\xi:\xi<\lambda\rangle$ is an ideal independent
sequence in $\bB^{\bf b}$. We may assume that $a_\xi=\bigwedge\limits_{k<
k_\xi} x^{t(\xi,k)}_{\alpha(\xi,k)}$ and $\alpha(\xi,k)\neq\alpha(\xi,k')$
whenever $k<k'<k_\xi$ (remember \ref{xxx}(2)). Also we may assume that
$k_\xi=k^*$ for all $\xi<\lambda$ (as $\cf(\lambda)>\omega$).

Fix $i<\cf(\lambda)$ for a moment.

After possibly re-enumerating the sequences $\langle\alpha(\xi,k):k<k^*
\rangle$, we may find a set $S_i\subseteq [\chi_i,\chi_{i+1})$, an ordinal
$\vare_i^*<\mu$, a sequence $\langle\nu^i_k:k<k^*\rangle$ of pairwise
distinct elements of $\mu^{\textstyle\vare^*_i}$, and $t^i_k\in\{0,1\}$ and
$j^i_k<\cf(\lambda)$ (for $k<k^*$) such that 
\begin{enumerate}
\item[(i)]   $S_i$ is unbounded in $\chi_{i+1}$, 
\item[(ii)]  $t(\xi,k)=t^i_k$ and $j(\alpha(\xi,k))=j^i_k$ for all $\xi\in
S_i$ and $k<k^*$,
\item[(iii)] $\nu^i_k\vartriangleleft\eta_{\alpha(\xi,k)}$ for $k<k^*$ and 
$\xi\in S_i$, 
\item[(iv)]  $\langle\langle\alpha(\xi,k):k<k^*\rangle:\xi\in S_i\rangle$
is a $\Delta$--system of sequences with heart $\langle\alpha^i_k:k<k(i)
\rangle$, 
\item[(v)]   the sequence $\langle\alpha(\xi,k):\xi\in S_i\rangle$ is
strictly increasing for $k(i)\leq k<k^*$,  
\item[(vi)]  $j^i_k\geq i$ for $k(i)\leq k<k^*$ (it follows from
(ii)+(iv)).
\end{enumerate}
Next pick a set $S\subseteq [\cf(\lambda)]^{\textstyle\cf(\lambda)}$ such
that (possibly after some re-enumerations)
\begin{enumerate}
\item[(vii)]  $k(i)=k^+$, $t^i_k=t_k$, $\vare^*_i=\vare^*$ and $\nu^i_k= 
\nu^*_k$ for  $k<k^*$, $i\in S$,
\item[(viii)] $\langle\langle\alpha^i_k:k<k^+\rangle:i\in S\rangle$ is a
$\Delta$--system of sequences with heart $\langle\alpha_k:k<k^{**}\rangle$, 
\item[(ix)]   $\langle\langle j^i_k: k<k^*\rangle: i\in S\rangle$ is a
$\Delta$--system of sequences with heart $\langle j_k:k\in w\rangle$,
$w\subseteq k^*$.
\end{enumerate}
Note that then $k^{**}\subseteq w\subseteq k^+$. Also, possibly further
shrinking $S$ and the $S_i$'s (for $i\in S$), we may demand that 
\begin{enumerate}
\item[(x)]   if $i_1<i_2$, $i_1,i_2\in S$, then $j^{i_1}_k<i_2$ (for
$k<k^*$),
\item[(xi)]  if $i_1,i_2\in S$ are distinct, $\xi_1\in S_{i_1}$ and $\xi_2
\in S_{i_2}$, then
\[\{\alpha(\xi_1,k):k<k^*\}\cap \{\alpha(\xi_2,k):k<k^*\}=\{\alpha_k:k<
k^{**}\}.\]
\end{enumerate}
Let $S^*=\bigcup\limits_{i\in S} S_i$. For $\vare<\mu$ and $k^+\leq k<k^*$
let
\[\begin{array}{r}
S^L_{\vare,k}=\big\{\xi\in S^*:(\forall\zeta\in S^*)\big(\vare>\lh(
\eta_{\alpha(\xi,k)}\cap\eta_{\alpha(\zeta,k)})\ \mbox{ or }\
\eta_{\alpha(\xi,k)}\cap\eta_{\alpha(\zeta,k)}\notin A\ \mbox{ or }\ \\
\eta_{\alpha(\xi,k)}\leq_\lex \eta_{\alpha(\zeta,k)}\big)\big\},\\
S^R_{\vare,k}=\big\{\xi\in S^*:(\forall\zeta\in S^*)\big(\vare>\lh(
\eta_{\alpha(\xi,k)}\cap\eta_{\alpha(\zeta,k)})\ \mbox{ or }\
\eta_{\alpha(\xi,k)}\cap\eta_{\alpha(\zeta,k)}\notin A\ \mbox{ or }\ \\
\eta_{\alpha(\zeta,k)}\leq_\lex \eta_{\alpha(\xi,k)}\big)\big\}.
  \end{array}\]
We claim that both $|S^L_{\vare,k}|<\lambda$ and $|S^R_{\vare,k}|<
\lambda$. Why? Assume, e.g., $|S^L_{\vare,k}|=\lambda$. Note that, by
(v)+(vi)+(x), $\alpha(\xi,k)<\alpha(\zeta,k)$ for $\xi<\zeta$ from $S^*$. 
Pick $\nu\in\mu^{\textstyle\vare}$ and a set $X\in [S^L_{\vare,k}]^{
\textstyle\lambda}$ such that $(\forall\xi\in X)(\nu\vartriangleleft
\eta_{\alpha(\xi,k)})$. By \ref{Def3.2}(1)(c), there are distinct $\xi,\zeta
\in X$ such that $\eta_{\alpha(\xi,k)}\cap\eta_{\alpha(\zeta,k)}\in
A$. Clearly $\lh(\eta_{\alpha(\xi,k)}\cap\eta_{\alpha(\zeta,k)})\geq\vare$
and we easily get a contradiction with $\xi,\zeta\in
S^L_{\vare,k}$. Similarly for $S^R_{\vare,k}$.

For $k^+\leq k<k^*$ let
\[\begin{array}{ll}
S^\otimes_k=\Big\{\xi\in S^*:&\mbox{for all }\vare<\mu\mbox{ there exists } 
\zeta\in S^*\mbox{ such that }\eta_{\alpha(\xi,k)}<_\lex\eta_{\alpha(\zeta,
k)},\\
&\mbox{and }\vare\leq\lh(\eta_{\alpha(\xi,k)}\cap\eta_{\alpha(\zeta,k)})\
\mbox{ and }\ \eta_{\alpha(\xi,k)}\cap\eta_{\alpha(\zeta,k)}\in A,\\
&\mbox{and }\\
&\mbox{for all }\vare<\mu\mbox{ there exists }\zeta\in S^*\mbox{ such that
}\eta_{\alpha(\zeta,k)}<_\lex\eta_{\alpha(\xi,k)},\\
&\mbox{and }\vare\leq\lh(\eta_{\alpha(\xi,k)}\cap\eta_{\alpha(\zeta,k)})\
\mbox{ and }\ \eta_{\alpha(\xi,k)}\cap\eta_{\alpha(\zeta,k)}\in A\Big\}. 
  \end{array}\]
Note that $S^*\setminus S^\otimes_k=\bigcup\limits_{\vare<\mu} (S^L_{\vare,
k}\cup S^R_{\vare,k})$, and hence $|S^*\setminus S^\otimes_k|<\lambda$ for
each $k\in [k^+,k^*)$. 

Fix distinct $\xi^*,\xi_*\in\bigcap\limits_{k=k^+}^{m-1} S^\otimes_k$ such
that $j(\xi^*)=j(\xi_*)$. For each $k\in [k^+,k^*)$ pick $\xi^k_1,\xi^k_2,
\xi^k_3\in S^*\setminus\{\xi^*,\xi_*\}$ such that 
\[\begin{array}{c}
\nu^*_k\vartriangleleft\eta_{\alpha(\xi^*,k)}\cap\eta_{\alpha(\xi^k_1,k)}
\vartriangleleft\eta_{\alpha(\xi^*,k)}\cap\eta_{\alpha(\xi^k_2,k)}
\vartriangleleft\eta_{\alpha(\xi^*,k)}\cap\eta_{\alpha(\xi^k_3,k)},\\
\eta_{\alpha(\xi^*,k)}\cap\eta_{\alpha(\xi^k_1,k)},
\eta_{\alpha(\xi^*,k)}\cap\eta_{\alpha(\xi^k_2,k)}\in A,\\
\eta_{\alpha(\xi^k_1,k)}<_\lex\eta_{\alpha(\xi^*,k)}<_\lex\eta_{\alpha(
\xi^k_2,k)}.
  \end{array}\]
Now look: letting $\alpha_k=\alpha(\xi^*,k)$, $\{\alpha_{\ell,k}:\ell<\ell^* 
\}$ be the suitable enumeration of $\big\{\alpha(\xi^{k'}_n,k):k^+\leq
k'<k^*\ \&\ n\in \{1,2,3\}\big\}\cup\{\alpha(\xi_*,k)\}$, and $\sigma_k=
\nu^*_k$, we get that the clauses $(\alpha)$--$(\gamma)$ of \ref{clx2} are
satisfied. Hence  
\[a_{\xi^*}=\bigwedge_{k<k^*} x^{t_k}_{\alpha(\xi^*,k)}\leq \bigwedge_{k<
k^*} x^{t_k}_{\alpha(\xi_*,k)}\vee\bigvee_{n=1}^3\bigvee_{k'=k^+}^{k^*-1}
\bigwedge_{k<k^*} x^{t_k}_{\alpha(\xi^{k'}_n,k)}=a_{\xi_*}\vee
\bigvee_{n=1}^3\bigvee_{k'=k^+}^{k^*-1}a_{\xi^{k'}_n}.\]
Since clearly $\xi^*\notin\{\xi_*\}\cup\{\xi^{k'}_n:k^+\leq k'<k^*, n=1,2,
3\}$, we get a contradiction. 
\end{proof}

\begin{claim}
\label{clx4}
If ${\bf b}$ is a base$^+$ then also $\hL_{(7)}^+(\bB^{\bf b})=\hd^+_{(5)}
(\bB^{\bf b})=\lambda$.
\end{claim}

\begin{proof}[Proof of the claim]
It is similar to \ref{clx3}. Suppose that $\langle a_\xi:\xi<\lambda\rangle$
is a right separated sequence in $\bB^{\bf b}$. Like before we may assume
that $a_\xi=\bigwedge\limits_{k<k^*} x^{t(\xi,k)}_{\alpha(\xi,k)}$ and
$\alpha(\xi,k)\neq\alpha(\xi,k')$ whenever $k<k'<k^*$. Next we apply the
same ``cleaning procedure'' as in \ref{clx3} getting $S,S_i,\vare^*,\nu^*_k,
t_k,j^i_k$ etc such that clauses (i)---(xi) are satisfied. We let $S^*=
\bigcup\limits_{i\in S} S_i$ and for $\vare<\mu$ and $k^+\leq k<k^*$ we
define  
\[\begin{array}{r}
S^+_{\vare,k}=\big\{\xi\in S^*:(\forall\zeta\!\in\!S^*\cap\xi)\big(\vare\!
>\!\lh(\eta_{\alpha(\xi,k)}\cap\eta_{\alpha(\zeta,k)})\ \mbox{ or }\
\eta_{\alpha(\xi,k)}\cap\eta_{\alpha(\zeta,k)}\notin A\quad\\ 
\mbox{or }\eta_{\alpha(\xi,k)}\leq_\lex\eta_{\alpha(\zeta,k)}\big)\big\},\\ 
S^-_{\vare,k}=\big\{\xi\in S^*:(\forall\zeta\!\in\!S^*\cap\xi)\big(\vare\!
>\!\lh(\eta_{\alpha(\xi,k)}\cap\eta_{\alpha(\zeta,k)})\ \mbox{ or }\
\eta_{\alpha(\xi,k)}\cap\eta_{\alpha(\zeta,k)}\notin A\quad\\ 
\mbox{or }\eta_{\alpha(\zeta,k)}\leq_\lex\eta_{\alpha(\xi,k)}\big)\big\}.
  \end{array}\]
Then both $|S^+_{\vare,k}|<\lambda$ and $|S^-_{\vare,k}|<\lambda$. [It is
like before: assume, e.g., $|S^+_{\vare,k}|=\lambda$. Pick $\nu\in
\mu^{\textstyle\vare}$ and a set $X\in [S^+_{\vare,k}]^{ \textstyle\lambda}$
such that $(\forall\xi\in X)(\nu\vartriangleleft\eta_{\alpha(\xi,k)})$. Note
that $\alpha(\zeta,k)<\alpha(\xi,k)$ for $\zeta<\xi$ from $S^*$. Use
\ref{Def3.2}(2)(c$^+$) to find $\zeta<\xi$, both from $X$, such that 
$\eta_{\alpha(\zeta,k)}\cap\eta_{\alpha(\xi,k)}\in A$ and $\eta_{\alpha(
\zeta,k)}<_\lex\eta_{\alpha(\xi,k)}$. A clear contradiction.]

Next for $k^+\leq k<k^*$ we let $S^\otimes_k=S^*\setminus\bigcup\limits_{
\vare<\mu}(S^+_{\vare,k}\cup S^-_{\vare,k})$. Choose $\xi_*<\xi^*$ from
$\bigcap\limits_{k=k^+}^{m-1} S^\otimes_k$ such that $j(\xi^*)=j(\xi_*)$.
And next for each $k\in [k^+,k^*)$ pick $\xi^k_1,\xi^k_2,\xi^k_3\in S^*\cap
\xi^*$ like those in the proof of \ref{clx3}. Finish in the same way.
\end{proof}
\end{proof}

\shlhetal
\end{document}